\documentclass{elsarticle}

\usepackage{lineno,hyperref}
\modulolinenumbers[1]

\usepackage{amsopn}	% \DeclareMathOperator

\usepackage{cleveref}
\usepackage{amssymb}
\usepackage{amsmath}
\usepackage{amsthm}
\usepackage{subcaption}
\usepackage{marginnote}
\usepackage{algorithm}
\usepackage{algorithmic}
\usepackage{multirow}
\usepackage{array}

\usepackage{graphicx}
\usepackage{color}
\usepackage{calc}

\graphicspath{ {./images/} }

\newcommand{\E}[2][]{\mathbb{E}_{{#1}} \left[ {#2} \right] }
\newcommand{\Var}[2][]{\mathbb{V}ar_{{#1}} \left[ {#2} \right] }
\newcommand{\Cov}[2][]{\mathbb{C}ov_{{#1}} \left[ {#2} \right] }
\newcommand{\Prob}[2][]{\mathbb{P}_{{#1}} \left[ {#2} \right] }
\newcommand{\norm}[2]{ \left\| #1 \right\|_{#2} }
\newcommand{\change}[1]{{\color{black}#1}}
\newcommand{\cchange}[1]{{\color{black}#1}}
\newcommand{\ccchange}[1]{{\color{black}#1}}

\newcolumntype{C}[1]{>{\centering\let\newline\\\arraybackslash\hspace{0pt}}m{#1}}

\DeclareMathOperator{\Tr}{Tr}
\DeclareMathOperator*{\argmin}{arg\,min}

\newtheorem{Definition}{Definition}
\newtheorem{Lemma}{Lemma}

\everymath{\displaystyle}
\journal{arXiv}

\begin{document}

\begin{frontmatter}

\title{Slow-scale split-step tau-leap method for stiff stochastic chemical systems}


\author[mtsu_address,ornl_address]{Viktor Reshniak\corref{mycorrespondingauthor}}
\cortext[mycorrespondingauthor]{Corresponding author}
\ead{reshniakv@ornl.gov}

\author[mtsu_address]{Abdul Khaliq}
\ead{Abdul.Khaliq@mtsu.edu}

\author[wiu_address]{David Voss}
\ead{d-voss1@wiu.edu}

\address[mtsu_address]{Department of Mathematical Sciences and Center for Computational Science, Middle Tennessee State University, Murfreesboro, TN 37132, USA}
\address[ornl_address]{Computer Science and Mathematics Division, Oak Ridge National Laboratory, Oak Ridge, TN 37831, USA}
\address[wiu_address]{Professor Emeritus, Department of Mathematics, Western Illinois University, Macomb, IL 61455, USA}

\begin{abstract}
Tau-leaping is a family of algorithms for the approximate simulation of the discrete state continuous time Markov chains.
A motivation for the development of such methods can be found, for instance, in the fields of chemical kinetics and systems biology.
It is known that the dynamical behavior of biochemical systems is often intrinsically stiff representing a serious challenge for their numerical approximation.
The naive extension of stiff deterministic solvers to stochastic integration often yields numerical solutions with either impractically large relaxation times or incorrectly resolved covariance.
In this paper, we propose a splitting heuristic which helps to resolve some of these issues.
\change{The proposed integrator contains a number of unknown parameters which are estimated for each particular problem from the moment equations of the corresponding linearized system.}
We show that this method is able to reproduce the exact mean and variance of the linear scalar test equation and demonstrates a good accuracy for the arbitrarily stiff systems at least in the linear case.
The numerical examples for both linear and nonlinear systems are also provided, and the obtained results confirm the efficiency of the considered splitting approach.
\end{abstract}


\end{frontmatter}

% \linenumbers
%\newpage

\section{Introduction}
\label{sec:intro}

\change{
The importance of stochastic fluctuations in biochemical processes has been established and confirmed by a large score of experimental and theoretical studies \cite{McAdams1999,Raj2008,Symmons2016}.
It is known that at the cellular level, even a single molecule can drastically impact the outcome of reactions.
Markov chain models can account for both discrete and stochastic nature of such processes and provide convenient tools for the mathematical description of stochastic chemical kinetics \cite{Erdi1989, Gillespie1991, Van1992}.
For example, the evolution of probabilities in Markovian kinetic networks is entirely specified by the (chemical) master equation (CME).
Unfortunately, its dimension grows exponentially with the number of molecular species restricting the domain of its application to very small or simple systems \cite{Gillespie1991,Hegland2007,Jahnke2007}. % with known exact solutions
% Stochastic simulation techniques effectively address this issue but require a large number of repetitive simulations of all reactions in the network. 

\ccchange{Simulation techniques provide an alternative to the master equation by generating large ensembles of possible temporal state paths.
% The classical 
Gillespie stochastic simulation algorithm (SSA) \cite{Gillespie1977} and its efficient implementations like the next reaction method \cite{Gibson2000}, the optimized direct method \cite{Cao2004_b} or the sorting direct method \cite{McCollum2006} keep track of every reaction event and are essentially exact.
Approximate simulation algorithms allow to skip over multiple reactions in a single leaping step and are often cheaper than SSA. % at the cost of losing accuracy.
The original tau-leaping algorithm in \cite{Gillespie2001} uses Poisson random variables to estimate the number of multiple reactions which can occur during the step of size $\tau$.
It is based on the assumption that $\tau$ is small enough to guarantee that the intensities of the reactions do not change drastically during the leaping step.
Several tau selection strategies have been proposed in the literature to satisfy this leaping condition \cite{Petzold2003,Petzold2006} and to avoid negative populations \cite{Tian2004,Petzold2005,Chatterjee2005,Moraes2014_b}.

% It is common for biochemical systems to contain both fast and slow reactions.
Many biochemical systems contain both fast and slow reactions.
% Fast reactions occur much more frequently than the others forcing tiny leaping steps and significantly reducing computational attractiveness of the approximate simulation methods.
It is known that transitory and highly reactive species involved in fast reactions may reach equilibrium and be asymptotically at a steady state within the coarse time scale of slow reactions \cite{Goutsias2005,Rao2003}.
Sampling fast intermediate species from their quasi-stationary distributions using SSA is redundant on time intervals much exceeding the appropriate relaxation times. 
% SSA is exact but very expensive for such systems while the approximate tau-leaping methods require tiny steps in order to remain stable and accurate and are also inefficient.
% SSA is exact but extremely inefficient for such systems.
Tau-leaping methods applied to such systems usually require tiny steps in order to remain stable and are also inefficient.}
This phenomenon is known as stiffness and is a well-studied topic in the numerical analysis of deterministic dynamical systems.
However, unlike deterministic dynamics, fast stochastic fluctuations represent an essential qualitative feature of  biochemical systems and not an artifact of a numerical integrator.
Efficient stiff solvers must be able to skip over the fast and stable reactions while capturing their stochastic influence on the slow species.
% This leads to the impractically high computational cost of the simulation algorithms when some reactions occur much more frequently than the others.
% This leads to the impractically high computational cost of the simulation algorithms when the dynamics of interest happens on the slow scale of the system with the fastest scale being stable.

% This phenomenon is known as stiffness and is a well studied topic in the numerical analysis of deterministic dynamical systems. 
% However, unlike deterministic dynamics, fast stochastic fluctuations represent an important qualitative feature of  biochemical systems and not an artifact of a numerical integrator.
% Efficient stiff solvers must be able to skip over the fast reactions while capturing their stochastic influence on the slow species.

% It is known that transitory and highly reactive species often reach equilibrium and are asymptotically at steady state within the coarse time scale of slow reactions \cite{Goutsias2005,Rao2003}.
% Using implicit solvers is the standard way of coping with stiffness in deterministic equations.
\ccchange{Stiffness in deterministic equations is efficiently handled with implicit solvers.}
The naive application of implicit integrators to stochastic systems yields numerical solutions with either unfeasibly large relaxation times or severely overdamped variance \cite{Rathinam2003, Cao2004}. %due to the overdamping property of implicit methods .
% This happens due to the overdamping property of implicit methods and their failure to correctly resolve invariant distributions of the fast variables.
Several approaches to cope with this issue have been proposed in the literature.}
Quasi-steady state or partial equilibrium approximation employ the splitting of reactions into fast and slow groups which are then integrated separately on their corresponding time scales \cite{Goutsias2005, Cao2005_3, Rao2003, Mastny2007}.
% Quasi-steady state or partial equilibrium approximation was employed for the dimension reduction of the singularly perturbed reaction networks in \cite{Goutsias2005, Cao2005_3, Rao2003, Mastny2007}. %, Han2014}.
% % It is based on separate integration of the fast and slow reactions on their corresponding time scales.
% It is based on the splitting of reactions into fast and slow groups which are then integrated separately on their corresponding time scales.
% These ideas were implemented in the so-called slow-scale and nested algorithms \cite{Cao2005, Cao2008, Liu2007}.
Slow-scale and nested algorithms implement these ideas in \cite{Cao2005, Cao2008, Liu2007}.
Other methods which utilize the splitting of reaction channels were considered in \cite{Hu2012, Jahnke2010, Engblom2015, E2005, Haseltine2002, Samant2007}.
Approaches which do not require explicit separation of scales have also been proposed.
For example, interlacing \cite{Rathinam2003,Cipcigan2011} and projective \cite{Lu2012} strategies were used to restore the overly damped stochastic fluctuations by interchanging large implicit steps with short explicit bursts.

Split-step methods have also been proposed for the numerical integration of stiff stochastic systems.
% Split-step methods form another class of integrators which attack the multi-scale nature of stochastic systems indirectly.
% The splitting of the original system into different parts provides the flexibility of choosing a combination of subsystem solvers with the most desirable set of properties.
For example, in \cite{Yang2011}, the splitting was used in conjunction with the Minkowski-Weyl decomposition to determine the joint distribution of the reaction count vector satisfying the nonnegativity and integrality conditions on the updated state.
Similarly, in \cite{Yang2013}, the joint distribution was approximated by a conditional Gaussian with the mean and covariance evaluated using the local central limit approximation. 
Additionally, split-step methods have been studied within the framework of stochastic Runge-Kutta integrators.
Chebyshev S-ROCK scheme for the systems with discrete noise was introduced in \cite{Abdulle2010}; it was shown to have excellent stability in the mean but failed to resolve the variance.
% It was shown that the method has excellent stability in the mean but the variance behaves poorly.
It was pointed out in \cite{Rue2010} that the accuracy of this scheme can be considerably improved by the appropriate choice of Runge-Kutta coefficients satisfying certain optimality conditions.

\change{
In this paper, we propose the new implicit split-step tau-leaping method which is both stable \ccchange{in the mean} and can accurately resolve stationary distributions of the \ccchange{chemical species involved in fast reactions} at least in the linear case.}
The single step of the proposed scheme has the form

\begin{align*}%\label{eq:composition}
	Y_{n+1} = \Xi(\tau)(Y_n),
\end{align*}
where the increment function 

\begin{align*}%\label{eq:increment}
	\Xi(\tau) = \Phi_2( (1-\boldsymbol{\theta}) \tau, \boldsymbol{\eta}_2 ) \circ S(\tau) \circ \Phi_1( \boldsymbol{\theta} \tau, \boldsymbol{\eta}_1 )
\end{align*}
consists of two implicit deterministic steps $\Phi_1$, $\Phi_2$ and the single tau-leaping step~$S$.
The functions $\Phi_1$ and $\Phi_2$ can be very general but we limit our attention to the classical theta scheme.
\change{The main result of the paper is the choice of the parameters $\boldsymbol{\eta}_1$, $\boldsymbol{\eta}_2$ and $\boldsymbol{\theta}$ which admits the accurate simulation of the fast and stable reactions on the coarse time grids.}
Our approach is conceptually similar to the one in \cite{Rue2010} as we also select the optimal parameters in a problem dependent manner.
The form of the splitting and the proposed parameter selection strategy, however, are different and, to the best of our knowledge, have not been applied to the equations of stochastic chemical kinetics before.

\ccchange{
% It is worth noting that in order to perform statistical inference for the quantities of interest associated with the state of the chemical system, stochastic simulation techniques are usually combined with Monte Carlo methods.
It is worth noting that stochastic simulation algorithms are usually combined with Monte Carlo methods.
A popular tool for reducing the computational burden of Monte Carlo sampling is provided by the multilevel techniques which rely on the availability of the sequence of approximations with increasing fidelity \cite{Giles2015}.
This approach was recently extended to the discrete state Markov chain models in \cite{Anderson2012,Moraes2015,Lester2016,BenHammouda2017}.
The proposed split-step method can potentially be used as a low-fidelity integrator at coarse levels of the discretization hierarchy improving the overall cost of the multilevel estimators.
% The proposed split-step method can be used as a low-fidelity integrator increasing the number of accessible coarse levels of the discretization hierarchy and improving the overall cost of the multilevel estimators.
% Multilevel Monte Carlo method is a popular tool for reducing the computational burden of random sampling by considering the sequence of approximations with increasing fidelity and .
% Multilevel estimators is a powerful tool for reducing the computational cost of the statistical estimators  rely on the availability of the sequence of approximations with increasing fidelity. 
}

The paper is organized as follows.
In \cref{sec:background}, we give an overview of the equations of stochastic chemical kinetics. %, discuss the sources of stiffness in these equations and the associated issues in their numerical treatment.
In \cref{sec:methods}, we review the classical simulation algorithms and introduce the new splitting heuristic.
\Cref{sec:stability} contains the main result of the paper; we conduct the comparative stability analysis of the classical and the proposed schemes and give the algorithm for estimating the parameters.
Finally, the numerical examples for both linear and nonlinear problems are given in \cref{sec:num_examples}.

\section{Stochastic chemical kinetics}
\label{sec:background}

Consider a thermally equilibrated chemical system which consists of $N$ well-stirred molecular species $\{ S_i, i=1,...,N \}$ interacting through $R$ reactions channels $\{ \mathcal{R}_r, r=1,...,R \}$

\begin{align*}
	\mathcal{R}_1 : \quad &\nu_{11}^{-} S_1 + \nu_{21}^{-} S_2 + ... +  \nu_{N1}^{-} S_N  \xrightarrow{c_1} \nu_{11}^{+} S_1 + \nu_{21}^{+} S_2 + ... +  \nu_{N1}^{+} S_N,
	\\
	\mathcal{R}_2 : \quad &\nu_{12}^{-} S_1 + \nu_{22}^{-} S_2 + ... +  \nu_{N2}^{-} S_N  \xrightarrow{c_2} \nu_{12}^{+} S_1 + \nu_{22}^{+} S_2 + ... +  \nu_{N2}^{+} S_N,
	\\
	\vdots
	\\
	\mathcal{R}_R : \quad &\nu_{1R}^{-} S_1 + \nu_{2R}^{-} S_2 + ... +  \nu_{NR}^{-} S_N  \xrightarrow{c_R} \nu_{1R}^{+} S_1 + \nu_{2R}^{+} S_2 + ... +  \nu_{NR}^{+} S_N.
\end{align*}
At any time instance, this system can be in precisely one of the states $X(t)=(X_1, X_2, ..., X_N)^T$, with $X_i$ denoting the number of molecules of type $S_i$. 
\ccchange{For chemical species with small to moderate molecular populations, intrinsic random fluctuations do not just average away and can lead to significant relative differences in the outcome of reactions even for initially identical system configurations. 
Evolution of the state vector for such systems can be modeled as a stochastic jump process \cite{Gillespie1977, Anderson2015}
}
% Under a reasonable assumption that collisions between molecules occur in a random manner, evolution of the state vector can be modeled as a stochastic jump process \cite{Gillespie1977, Anderson2015}

\begin{align}\label{eq:jump_proc}
	\begin{pmatrix}
		X_1(t) \\
		X_2(t) \\
		\vdots \\
		X_N(t) \\
	\end{pmatrix}
	=
	\begin{pmatrix}
		X_1(0) \\
		X_2(0) \\
		\vdots \\
		X_N(0) \\
	\end{pmatrix}
	+
	\begin{pmatrix}
		\nu_{11} & \nu_{12} & \hdots & \nu_{1R} \\
		\nu_{21} & \nu_{22} & \hdots & \nu_{2R} \\
		\vdots   &          & \ddots & \vdots   \\
		\nu_{N1} & \nu_{N2} & \hdots & \nu_{NR} \\
	\end{pmatrix}
	\begin{pmatrix}
		N_1(t) \\
		N_2(t) \\
		\vdots \\
		N_R(t) \\
	\end{pmatrix},
\end{align}
where each stoichiometric column vector $\nu_r = [\nu_{1r}^{+}-\nu_{1r}^{-},...,\nu_{Nr}^{+}-\nu_{Nr}^{-}]^T$ denotes the change in molecular populations from the reaction channel $\mathcal{R}_r$, and Markovian processes $N_r(t)$ count the number of corresponding reactions in the time interval $[0,t]$. 

Equation \eqref{eq:jump_proc} represents a continuous time Markov chain with the transition probabilities

\begin{align*}
	\Prob{ X(t+dt) - X(t) = \nu_r | X(t) } = a_r(X(t)) dt + o(dt), \qquad r = 1,..,R
\end{align*}
and with the associated chemical master equation

\begin{align}\label{eq:CME}
	\frac{\partial \mathbb{P}_{t,x}}{\partial t} &= -\sum_{r=1}^R a_r(x) \mathbb{P}_{t,x} + \sum_{r=1}^R a_r(x-\nu_r) \mathbb{P}_{t,x-\nu_r},
	\\ \nonumber
	\mathbb{P}_{0,x} &= p_x.
\end{align}
Equation \eqref{eq:CME} describes the change in time of the probability mass function $\mathbb{P}_{t,x} = \Prob{X(t)=x}$ of  every element $x \in \mathbb{S}$ from the state space $\mathbb{S} \subset \mathbb{Z}_{+}^N$ of the system.
According to the stochastic law of mass action, transition intensities (propensities) $a_r(x)$ are proportional to the reaction-rate constants $c_r$ and the number of distinct combinations of molecules of the source species in the corresponding reactions \cite{Anderson2015}

\begin{align*}
	a_r(x)= c_r \prod_{i=1}^N \binom{x_i}{\nu_{ri}^{-}}  = c_r \prod_{i=1}^N \frac{x_i!}{(x_i-\nu_{ri}^{-})! \nu_{ri}^{-}!}, \qquad r = 1,..,R.
\end{align*}
	
Although the master equation provides a complete probabilistic description of the chemical system at any time instance, it admits analytical solutions only in the exceptional simple cases and rarely can be solved numerically due to the exponential growth of its dimension with the number of species involved. 
In practice, however, it often suffices to know the evolution law of statistical moments \ccchange{since for any analytic function $g(x)$ of a random variable $x$, we have that

\begin{align*}
    \E{g(x)} = \sum_{n=0}^{\infty} \frac{g^{(n)}(0)}{n!} \E{x^n}.
\end{align*}
Additionally, multiscale stochastic simulation algorithms often require at most two first stationary moments of the fast species to compute the averaged slow dynamics \cite{Cao2005,Cao2005_3}.

}

Consider the expected value of the arbitrary function of the state 

\begin{align*}
	\E{f(X(t))} &= \sum_{x \in \mathbb{S}} f(x) \mathbb{P}_{t,x}.
\end{align*}
By differentiating the above expression and using \eqref{eq:CME}, we get

\begin{align}\label{eq:moment_evolution}
	\nonumber	
	\frac{\partial}{\partial t} \E{f(X(t))}
	&= \sum_{x \in \mathbb{S}} f(x) \frac{\partial \mathbb{P}_{t,x}}{\partial t} 
	= \sum_{r=1}^R \sum_{x \in \mathbb{S}} f(x) \Big( a_r(x-\nu_r) \mathbb{P}_{t,x-\nu_r} - a_r(x) \mathbb{P}_{t,x} \Big)
	\\ 
	&= \sum_{r=1}^R \sum_{x \in \mathbb{S}} \Big( f(x+\nu_r) - f(x) \Big) a_r(x) \mathbb{P}_{t,x}.
\end{align}
Note that the change of variables $(x-\nu_r) \to x$ in the first term in parentheses does not change the limits of summation because the sum is taken over all possible states and both $(x-\nu_r) \in \mathbb{S}$ and $x \in \mathbb{S}$.

In particular, we have

\begin{align*}
	&\frac{\partial}{\partial t} \E{X_i}
	= \sum_{r=1}^R \sum_{x \in \mathbb{S}} \Big( x_i + \nu_{i,r} - x_i \Big) a_r(x) \mathbb{P}_{t,x}
	= \sum_{r=1}^R \nu_{i,r} \E{ a_r(X) },
	\\
	&\frac{\partial}{\partial t} \Big( \E{X_i X_k} - \E{X_i} \E{X_k} \Big)
	\\	
	&= \sum_{r=1}^R \Big( 
	  \nu_{k,r} \big( \E{X_i a_r(X)} - \E{X_i} \big)
	+ \nu_{i,r} \big( \E{X_k a_r(X)} - \E{X_k} \big) 
	+ \nu_{i,r} \nu_{k,r} \E{a_r(X)} \Big).
\end{align*}
In matrix-vector notation, the above expressions take the form

\begin{align*}
	%\label{eq:mean_evol}
	\frac{\partial}{\partial t} \E{X(t)} &= \nu \cdot \E{ a\big(X(t)\big) },
	\\ %\label{eq:var_evol}
	\frac{\partial}{\partial t} \Cov{X(t),X(t)} &= \nu \cdot \Cov{a(X(t)),X(t)} + \Cov{X(t),a(X(t))} \cdot \nu^T 
	\\ \nonumber
	&+ \nu \cdot diag\left( \E{ a\big(X(t)\big) } \right) \cdot \nu^T.
\end{align*}
This system may not be closed due to possible higher moments of $X(t)$ appearing in the right-hand side of the equations.
\ccchange{However, the propensity functions of the zero and first order reactions can be combined into a single vector in the form of}
% However, for the particular case of the \ccchange{reactions up to first order}, the propensity functions are \ccchange{at most} linear

\begin{align}\label{eq:linear_prop}
	a(X) = C X + d,
\end{align}
where the nonzero elements $C_{rj}=c_r$ of the $R \times N$ matrix $C$ correspond to the rates of the conversion and degradation reactions $\mathcal{R}_r: S_j {\stackrel{c_{r}}{\rightarrow}} *$ and the $R$~-~dimensional vector $d$ contains the rates of the inflow reactions $\mathcal{R}_r: \emptyset {\stackrel{c_{r}}{\rightarrow}} *$. %, where $*$ denotes the arbitrary products of reactions.
In this case, the mean and covariance of the state vector $X(t)$ are defined by the following closed system

\begin{align}
	\label{eq:linear_mean_evol}
	\frac{\partial}{\partial t} \E{X(t)} &= \nu C \cdot \E{X(t)} + \nu d,
	\\ \label{eq:linear_var_evol}
	\frac{\partial}{\partial t} \Cov{X(t)} &= \nu C \cdot \Cov{X(t)} + \Cov{X(t)} \cdot (\nu C)^T 
	\\ \nonumber
	&+ \nu \cdot diag\left( C \cdot \E{X(t)} + d \right) \cdot \nu^T.
\end{align}

% In \eqref{eq:linear_prop}-\eqref{eq:linear_var_evol}, the nonzero elements $C_{rj}=c_r$ of the $R \times N$ matrix $C$ correspond to the rates of the conversion and degradation reactions $\mathcal{R}_r: S_j {\stackrel{c_{r}}{\rightarrow}} *$ and the $R$~-~dimensional vector $d$ contains the rates of the inflow reactions $\mathcal{R}_r: \emptyset {\stackrel{c_{r}}{\rightarrow}} *$, where $*$ denotes the arbitrary products of reactions.

\section{Simulation algorithms}
\label{sec:simulation}

\change{
As an alternative to the CME, the statistical quantities of interest associated with the state of the chemical system can be extracted directly from \eqref{eq:jump_proc} by simulating large ensembles of individual state paths.
In this section, we overview the classical algorithms for the exact and approximate path simulation and propose a new one.
}

\change{\subsection{Exact simulation. Stochastic simulation algorithm}}
\label{sec:exact_simulation}

The random time change formula allows reformulating the counting processes $N_r(t)$ in \eqref{eq:jump_proc} in terms of the unit rate Poisson processes $\mathcal{P}(t)$ \cite[Theorem~1.10]{Anderson2015}.
This results in the following stochastic integral equation

\begin{align}\label{eq:jump_Poiss}
	X(t) = X(0) + \sum_{r=1}^R \nu_r \mathcal{P}\left(\int_0^t a_r(X(s)) ds \right).
\end{align}
Taking into account that the holding times of the Poisson jump process are exponentially distributed random variables, the individual realizations of $X(t)$ can be generated with the stochastic simulation algorithm (SSA) as follows \cite{Gillespie1977}

\begin{align*}
	X(t_{n}+\tau) &= X(t_{n}) + \nu_r.
\end{align*}
The holding time and the index of the next reaction are calculated as

\begin{align*} %\label{eq:SSA}
	\tau &= -\frac{1}{\sum_{i=1}^{R} a_i(X(t_n))} \ln r_1,
	\\[0.5em]
	r &= \min \left\{ r:\sum_{i=1}^{r} a_i(X(t_n)) > r_2 \sum_{i=1}^{R} a_i(X(t_n)) \right\},
\end{align*}
where $r_1$, $r_2$ are the standard uniform random numbers.
There also exist several efficient implementations of this algorithm such as the next reaction method \cite{Gibson2000}, the optimized direct method \cite{Cao2004_b} and the sorting direct method \cite{McCollum2006}.

\ccchange{As was pointed out above, SSA is exact but redundant for sampling random state paths from the stationary distributions.}
Approximate path simulation algorithms can be more suitable in this case since they allow to skip over multiple reactions in a single step.

% \change{
% In chemical systems with time scale separation, transitory and highly reactive species often reach equilibrium and are asymptotically at steady state within the coarse time scale \cite{Goutsias2005, Rao2003}.
% It is obvious that the sampling of intermediate species from their quasi-stationary distributions is redundant on the time intervals which are larger than the appropriate relaxation times. 
% Stochastic simulation algorithm is exact but extremely inefficient for such systems.
% Approximate path simulation algorithms can be more suitable in this case since they allow to skip over multiple reactions in a single large step.}

\change{\subsection{Approximate simulation. Tau-leaping}}
\label{sec:methods}

Consider the integral equation \eqref{eq:jump_Poiss} over the time interval $t~\in~[t_n, t_{n+1}]$
\begin{align}\label{eq:martingale_noise}
	X(t_{n+1}) 
	&= X(t_n) + \sum_{r=1}^R \nu_r \mathcal{P}\left(\int\displaylimits_{t_n}^{t_{n+1}} a_r(X(s)) ds \right)
	\\ \nonumber
	&= X(t_n) + \sum_{r=1}^R \nu_r \int\displaylimits_{t_n}^{t_{n+1}} a_r(X(s)) ds + \sum_{r=1}^R \nu_r \overline{\mathcal{P}}\left(\int\displaylimits_{t_n}^{t_{n+1}} a_r(X(s)) ds \right)
\end{align}
with the driving martingale stochastic processes

$$ \overline{\mathcal{P}}\left(\int\displaylimits_{t_n}^{t_{n+1}} a_r(X(s)) ds \right)~=~\mathcal{P}\left(\int\displaylimits_{t_n}^{t_{n+1}} a_r(X(s)) ds \right)~-~\int\displaylimits_{t_n}^{t_{n+1}} a_r(X(s)) ds.$$

% The major difficulty in the construction of accurate numerical solutions to the above equation is associated with the approximation of the state dependent Poisson processes.
Approximation of the state dependent Poisson processes is the major difficulty in the construction of accurate numerical solutions to the above equation.
For slow reactions, it is reasonable to assume that the number of molecules and hence propensities do not change fast inside the interval of integration and one can use an explicit approximation of the integrals

$$ \overline{\mathcal{P}}\left(\int\displaylimits_{t_n}^{t_{n+1}} a_r(X(s)) ds \right) \approx \overline{\mathcal{P}}\Big( a_r(X(t_n)) \tau \Big), \qquad \tau=t_{n+1}-t_n.$$

Approximating the deterministic part of \eqref{eq:martingale_noise} with the theta scheme, one gets the family of drift implicit tau-leaping methods

\begin{align}\label{eq:theta}
	Y_{n+1}
	&= Y_n + \sum_{r=1}^R \nu_r \Big( (1-\theta) a_r(Y_{n}) + \theta a_r(Y_{n+1}) \Big) \tau + \sum_{r=1}^R \nu_r \overline{\mathcal{P}}\Big(a_r(Y_{n}) \tau \Big) 
	\\ \nonumber
	&= Y_n + \theta \sum_{r=1}^R \nu_r a_r(Y_{n+1}) \tau + \sum_{r=1}^R \nu_r \Big[ \mathcal{P}\Big(a_r(Y_{n}) \tau \Big) - \theta a_r(Y_{n}) \tau \Big],
\end{align}
where $Y_n$ is the numerical approximation of $X(t_n)$.
The classical explicit \cite{Gillespie2001}, implicit \cite{Rathinam2003} and trapezoidal \cite{Cao2004} tau-leaping methods correspond to the parameter values $\theta=0$, $\theta=1$ and $\theta=\frac{1}{2}$ respectively.
The error analysis of these schemes reveals their first order of weak convergence \cite{Li2007, Anderson2011a, Rathinam2005}.

Theta methods are implicit in their deterministic part but use the explicit approximation of the driving stochastic processes. 
Split-step methods form another class of integrators which allow incorporating implicitness into the stochastic part as well.
The standard split-step scheme reads as

\begin{align}\label{eq:spli_step}
  \hat{Y}_{n}
  %= Y_n &+ \sum_{r=1}^R \nu_r a_r(\hat{Y}_{n}) \tau,
  = Y_n &+ \sum_{r=1}^R \nu_r \Big( (1-\theta) a_r(Y_{n}) + \theta a_r(\hat{Y}_{n}) \Big) \tau,
  \\ \nonumber
  Y_{n+1}
  = \hat{Y}_{n} &+ \sum_{r=1}^R \nu_r \overline{\mathcal{P}}\left( a_r(\hat{Y}_{n}) \tau \right)
  = {Y}_{n} + \sum_{r=1}^R \nu_r {\mathcal{P}}\left( a_r(\hat{Y}_{n}) \tau \right).
\end{align}

In this paper, we propose the two-stage split-step method which combines ideas of the theta and the split-step methods.
% \eqref{eq:theta} and \eqref{eq:spli_step}.
The proposed scheme has the form

\begin{align}\label{eq:mod_split_step}
    \nonumber
    \hat{Y}_{n}
    = Y_n &+ \sum_{r=1}^R \nu_r \left( (1-\eta_{1r}) a_r(Y_n) + \eta_{1r} a_r(\hat{Y}_{n}) \right) (1-\theta_r) \tau,
    \\
    \tilde{Y}_{n}
    = \hat{Y}_{n} &+ \sum_{r=1}^R \nu_r \overline{\mathcal{P}}\left( a_r(\hat{Y}_{n}) \tau \right),
    \\ \nonumber
    Y_{n+1}
    = \tilde{Y}_n &+ \sum_{r=1}^R \nu_r \left( (1-\eta_{2r}) a_r(\tilde{Y}_n) + \eta_{2r} a_r(Y_{n+1}) \right) \theta_r \tau
\end{align}	
with the parameters $\eta_{1r}, \eta_{2r}, \theta_r \in [a,b]$ for some bounded $a$, $b$.
% The optimal selection of these parameters is discussed in section \ref{sec:stability}.
We discuss the optimal selection of these parameters in section \ref{sec:stability}.

\textbf{Remark.} 
Numerical schemes for the stochastic biochemical systems must ensure the non-negativity and integrality of the generated states.
The integrality for the scheme in \eqref{eq:mod_split_step} can be guaranteed if we set

\begin{align*}
	Y_{n+1}^{int}
	= Y_n  
	+ \sum_{r=1}^R \nu_r \Bigg\lfloor 
 	& \left( (1-\eta_{1r}) a_r(Y_n) + \eta_{1r} a_r(\hat{Y}_{n}) \right) (1-\theta_r) \tau 
	+ \overline{\mathcal{P}}\left( a_r(\hat{Y}_{n}) \tau \right)
	\\
	+ & \left( (1-\eta_{2r}) a_r(\tilde{Y}_n) + \eta_{2r} a_r(Y_{n+1}) \right) \theta_r \tau \Bigg\rceil,
\end{align*}	
where $\lfloor \cdot \rceil$ denotes the nearest integer function.
\ccchange{This choice ensures stoichiometrically realizable states since every state change is of the form $k_1\nu_1+...+k_R\nu_R$ for some non-negative integers $k_1,...,k_R$.}
The nonnegativity of the states can be achieved by applying the bounding procedure proposed in \cite[Section 2.2]{Rathinam2007}\ccchange{; it also follows stoichiometry of the system and does not alter the first order consistency of the tau-leaping methods by updating events with probabilities of order $O(\tau^2)$ or higher \cite{Rathinam2005}.}

% \change{
% Tau leaping methods are generally much faster than the stochastic simulation algorithm since they allow to skip over the fast reactions by integrating \eqref{eq:martingale_noise} on the time scale of the slow reactions. 
% However, this can also lead to the significant reduction of accuracy when the influence of the fast scale on the slow species is not resolved correctly.
% In the next section, we show that the classical tau-leaping methods are not well suited for this purpose while the proposed scheme demonstrates good results.
% }

\section{Moment stability analysis}
\label{sec:stability}

\ccchange{As was mentioned above, fast and highly reactive species often converge to the stable state within the coarse time scale of slow reactions.
Here stability is understood as the existence of a stationary distribution of the fast species conditioned on the ``frozen" slow species.
However, distributions are hard to study.
Analysis of statistical moments is often more tractable, especially in the linear case.
It has been shown that the classical tau-leaping methods fail at an accurate simultaneous resolution of the stationary mean and variance \cite{Cao2004}.}
\change{
% As was mentioned above, stable quasi-stationary distributions of the fast molecular species account for their stochastic influence on the slow species.
% For the obvious reason, working with exact analytical distributions is almost never possible.
% The analysis of statistical moments is often more tractable, especially in the linear case.
In this section, we provide the comparative linear moment stability analysis of the theta tau-leaping method \eqref{eq:theta} and the split-step scheme \eqref{eq:mod_split_step} and show the superiority of the latter for integrating fast and stable chemical reactions.
%in application to stiff chemical systems.
}

We will need the following definitions
\change{
\begin{Definition}[Reaction rate equation]\label{def:reaction_rate}
Reaction rate equation (RRE) is a deterministic ODE model for the evolution of molecular concentrations 
	\begin{align*}
		\frac{dz(t)}{dt} = \sum_{r=1}^R \nu_r \overline{a}_r(z(t)),
	\end{align*}
    where $z=\frac{X}{V}$, $\overline{a_r}=\frac{a_r}{V}$ and $V$ is the volume of the system.
    This model gives a good approximation of the dynamical behavior of chemical systems in the fluid limit such that $\lim_{V\to\infty} \frac{X(0)}{V} = const $ \cite{Ethier2009,Gillespie1991}. 
\end{Definition}
}

\cchange{
\begin{Lemma}[Law of total expectation]\label{def:law_total_mean}
The unconditional expectation and covariance of the random vector $Y_{n+1}$ can be evaluated as
	\begin{align*}
		\E{Y_{n+1}}   &= \E{\E{Y_{n+1} | Y_n}},
		\\[0.5em]
		\Cov{Y_{n+1}} &= \Cov{\E{Y_{n+1} | Y_n}} + \E{\Cov{Y_{n+1} | Y_n}}.
	\end{align*}
\end{Lemma}
}

\subsection{Motivating example}

As a motivation for the choice of the splitting in \eqref{eq:mod_split_step}, consider the reversible isomerization reaction %which serves as analog of the deterministic Dahlquist's test equation \cite{Cao2004}

\begin{align}\label{eq:test}
	S_1 \underset{c_2}{\stackrel{c_1}{\rightleftharpoons}} S_2.
\end{align}
This reaction has the linear propensities $a_1 = c_1 X_1$, $a_2 = c_2 X_2$ with the stoichiometric vectors $\nu_1 = (-1, 1)^T$ and $\nu_2=(1, -1)^T$.
The stationary distribution of \eqref{eq:test} is the binomial distribution \cite{Cao2004,Jahnke2007}

$$\Prob{X^*=x} = \frac{x_T!}{x_T!(x_T-x)!}q^x(1-q)^{x_T-x}$$ with the mean and the variance 

\begin{align*}
	\E{X^*}&= \frac{c_2}{\lambda}x_T,
	\qquad
	\Var{X^*} = \frac{c_1 c_2}{\lambda^2}x_T,
\end{align*} 
where $X^{*} = X_1(\infty)$, $\lambda=c_1+c_2$ and $x_T$ denotes the fixed total number of molecules, i.e., $x_T=X_1+X_2$.

Let $Y_n$ denote the numerical approximation of $X_1(t_n)$ and consider the theta method \eqref{eq:theta} applied to the test system \eqref{eq:test}.
The following condition ensures the global stability of the mean and the variance of the numerical solution (see \ref{sec:theta_stab})

\begin{align}\label{eq:theta_stab}
    |P(\lambda\tau)| = \left| \frac{1 - \lambda (1-\theta)\tau}{1+\lambda\theta\tau} \right| < 1.
\end{align}
\ccchange{We call $P(\lambda\tau)$ the propagation coefficient.} 

Assuming \eqref{eq:theta_stab} and letting $n\to\infty$, we get the mean and the variance of the corresponding stationary distribution

\begin{align}\label{eq:theta_stat_dist}
	\E{Y_{\infty}} 
	&= \E{X^{*}},
	\\ \nonumber
	\Var{Y_{\infty}} 
	&= A(\lambda\tau) \Var{X^{*}},
\end{align}
\ccchange{where $A(\lambda\tau) = \frac{2}{2+\lambda(2 \theta-1)\tau}$ is the stationary variance amplifier.}

\begin{figure}[!t]
	\centering
	\includegraphics[width=0.55\textwidth]{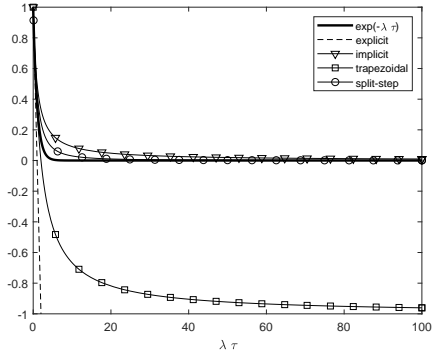} 
	\caption{Propagation coefficient of numerical schemes.}
	\label{fig:propagation_coef}
\end{figure}

Table \ref{tab:theta} shows the propagation coefficients and the stationary variance amplifiers of the theta scheme for different values of $\theta$.
It is seen that the stationary mean of the true solution is recovered correctly for all values of~$\theta$ when condition \eqref{eq:theta_stab} is satisfied. 
The stationary variance is resolved correctly only by the trapezoidal scheme with $\theta=\frac{1}{2}$, other values of $\theta$ result in either underdamped ($\theta<\frac{1}{2}$) or overdamped ($\theta>\frac{1}{2}$) variance.
Additionally, one can see \cchange{from Figure~\ref{fig:propagation_coef}} that the propagation coefficient of the trapezoidal scheme is close to one for large values of $\lambda \tau$ indicating slow relaxation rate of the numerical solution to the stationary state.
This limits the domain of application of the trapezoidal scheme to the systems with the time scale separation much larger than the relaxation time of the numerical solution.

\begin{table}[t]
\centering
\def\arraystretch{3}
\begin{tabular}{|c|p{\widthof{Propagation coefficients}}|p{\widthof{Variance amplifiers}}|c|}
\hline
 $\theta$ & \centering{Propagation coefficients \newline\newline $\frac{1-\lambda(1-\theta)\tau}{1+\lambda\theta\tau}$} & \centering{Variance amplifiers \newline\newline $\frac{2}{2+\lambda(2 \theta-1)\tau}$} & Stability condition \\ \hline \hline
 0 & \centering{$1-\lambda \tau $} & \centering{$\frac{2}{2 - \lambda \tau}$} & $\tau < \frac{\cchange{2}}{\lambda}$\\ 
 $\frac{1}{2}$ & \centering{$\frac{2-\lambda\tau}{2+\lambda\tau}$} & \centering{1} & unconditional\\ 
 1 & \centering{$\frac{1}{1+\lambda\tau}$} & \centering{$\frac{2}{2 + \lambda \tau}$} & unconditional \\ \hline
\end{tabular}
\caption{Propagation coefficients and stationary variance amplifiers of the theta method \eqref{eq:theta}. Parameter $\lambda$ is set to $\lambda=c_1+c_2$.}
\label{tab:theta}
\end{table}

Now consider the split-step method \eqref{eq:mod_split_step} applied to the system \eqref{eq:test}.
Let $\theta_r=\theta$, $\eta_{1r}=\eta_1$ and $\eta_{2r}=\eta_2$.
The global stability condition reads as (see \ref{sec:split_step_stab})
\begin{align*}
	|P(\lambda\tau)|
	=\left|  \left(\frac{1-(1-\eta_1)\lambda(1-\theta)\tau}{1+\eta_1\lambda(1-\theta)\tau} \right) \left( \frac{1-(1-\eta_2)\lambda\theta\tau}{1+\eta_2\lambda\theta\tau} \right)  \right| < 1.
\end{align*}	
Assuming this condition, the stationary moments of the numerical solution take the following values

\begin{align*}
	\E{Y(\infty)}
	&=  \E{X^{*}},
	\\
	\Var{Y(\infty)}
	&=  A(\lambda\tau) \Var{X^{*}}
\end{align*}
with 

$$
    A(\lambda\tau) = \frac{2\lambda\tau}{\displaystyle{\left( \frac{1+\eta_2\lambda\theta\tau}{1-(1-\eta_2)\lambda\theta\tau} \right)^2 - \left( \frac{1-(1-\eta_1)\lambda(1-\theta)\tau}{1+\eta_1\lambda(1-\theta)\tau} \right)^2}}.
$$

One can see that the split-step method also preserves the stationary mean of the true solution while the stationary variance depends on the choice of the parameters $\theta$, $\eta_1$, $\eta_2$. % and has more difficult behavior.
Table \ref{tab:split_step} provides the propagation coefficients and the variance amplifiers of the scheme \eqref{eq:mod_split_step} for different values of these parameters.
It is seen that the proposed scheme recovers the stability properties of the theta tau-leaping method when $\eta_1=0$, $\eta_2=1$ and of the standard split-step scheme \eqref{eq:spli_step} when $\theta=0$, $\eta_1=1$.
\cchange{By rewriting the variance amplifier as a rational function in $\lambda\tau$, it is also easy to check that the only choice of fixed parameters for which $\Var{Y(\infty)}=\Var{X^*}$ uniformly in $\lambda\tau$ is given by $\eta_1=0$, $\eta_2=1$, $\theta=\frac{1}{2}$ which corresponds to the trapezoidal theta scheme.}

\begin{table}[t]
\centering
\def\arraystretch{3}
\begin{tabular}{|c|C{0.24\textwidth}|C{0.34\textwidth}|C{0.18\textwidth}|}
% \begin{tabular}{|c|C{0.24\textwidth}|c|c|}
\hline
$[\theta,\eta_1,\eta_2]$ & Propagation \newline coefficients  & Variance \newline amplifiers & \cchange{Stability \newline condition} \\ \hline \hline
% $\theta$ & $\eta_1$ & $\eta_2$ & coefficients & amplifiers & \cchange{condition} \\ \hline \hline
$[\theta,1,1]$ & $\frac{1}{1+z+\theta(1-\theta)z^2}$ & $\frac{2 z}{(1+\theta z)^2 - \displaystyle{\frac{1}{(1+(1-\theta)z)^2}}}$ & \cchange{unconditional} \\
$[\theta,0,1]$ & $\frac{1-(1-\theta)z}{1+\theta z}$ & $\frac{2}{2+(2 \theta-1)z}$ & \cchange{$(1-2\theta) z < 2$} \\
$[0,\eta_1,\eta_2]$ & $\frac{1-(1-\eta_1)z}{1+\eta_1 z}$ & $\frac{2 (1+\eta_1z)^2}{2+z (1-2\eta_1+2\eta_1^2)}$ & \cchange{$(1-2\eta_1) z < 2$} \\
 \hline
\end{tabular}
\caption{Propagation coefficients and stationary variance amplifiers of the split-step method \eqref{eq:spli_step} with $z=\lambda\tau$ and $\lambda=c_1+c_2$.}
\label{tab:split_step}
\end{table}

\begin{figure}[!t]
	\centering
	\includegraphics[width=0.49\textwidth]{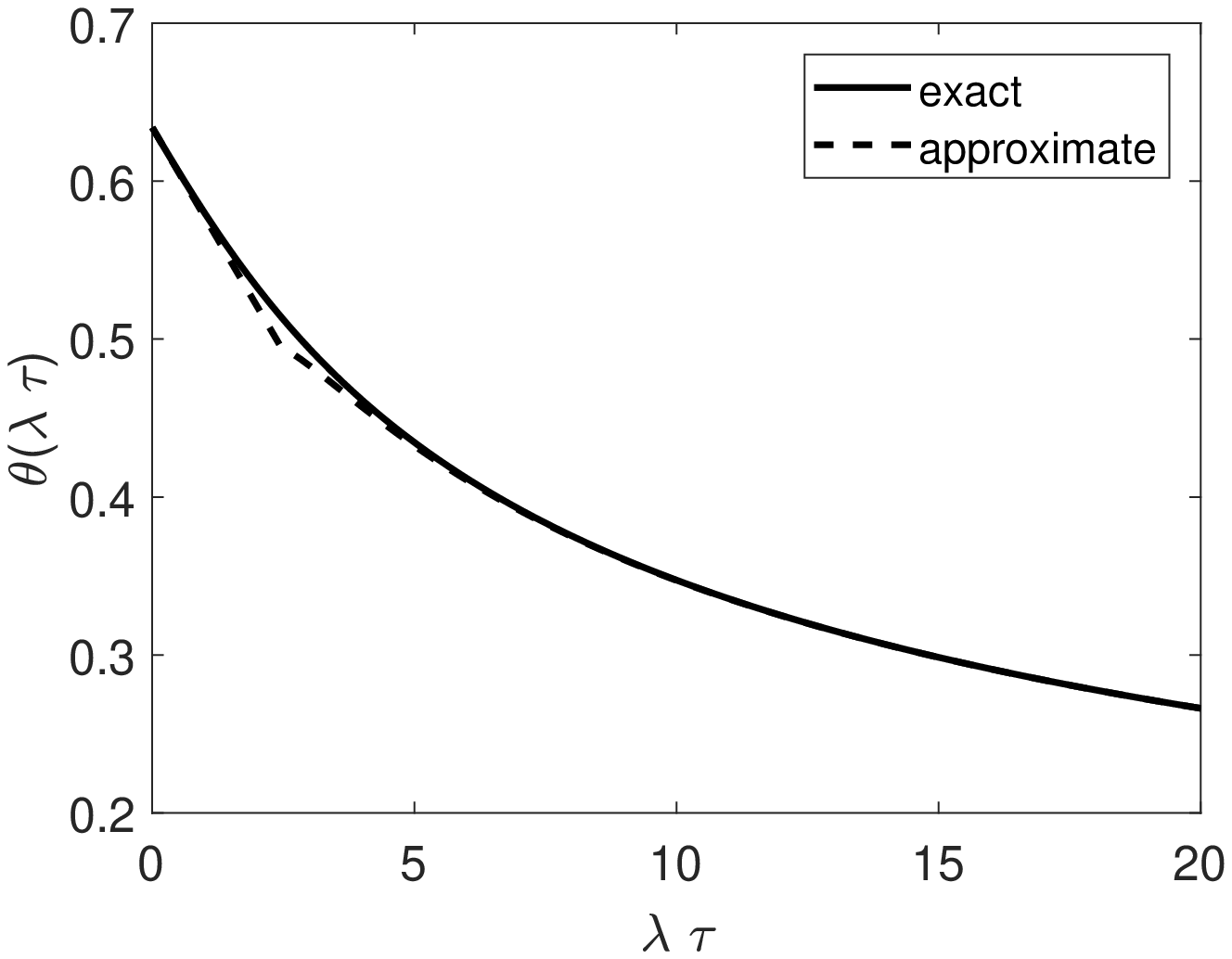} 
	\includegraphics[width=0.49\textwidth]{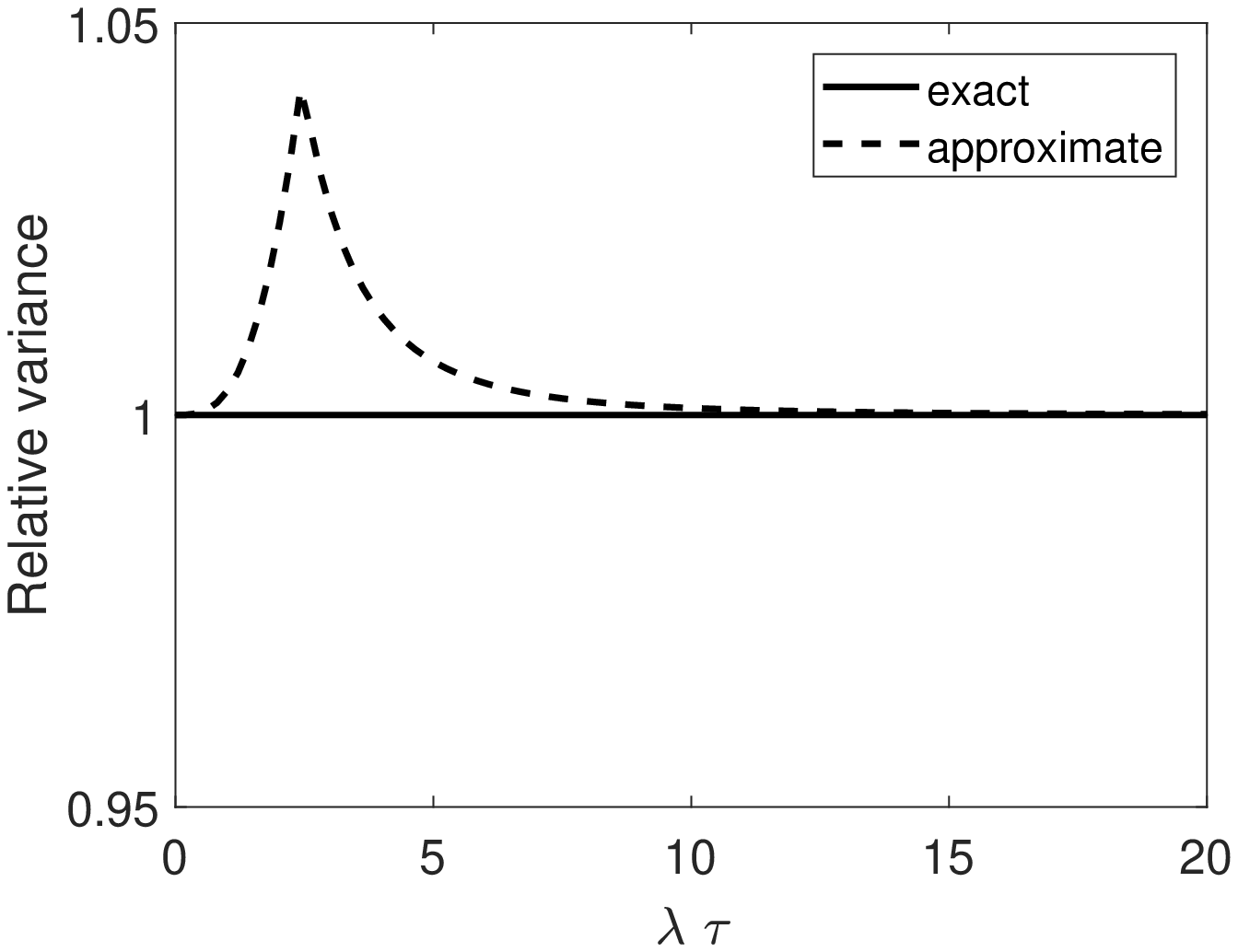}
	\caption{Parameter $\theta$ and the corresponding variance amplifier.}
	\label{fig:theta_choice}
\end{figure}

\cchange{It is worth noting that with $\eta_1=\eta_2=1$, the global stability condition for the split-step method is satisfied uniformly in $\theta$ and the propagation coefficient tends to zero for large values of $\lambda \tau$ producing numerical solutions with small relaxation time.
Figure \ref{fig:propagation_coef} illustrates the propagation coefficients from Tables \ref{tab:theta} and \ref{tab:split_step} for different numerical schemes along with the exponential decay $e^{-\lambda t}$ of the exact mean $\E{X_1(t)}$.
It is clear that the split-step method with $\eta_1=\eta_2=1$ gives the best approximation of the true exponential decay when compared to other considered schemes.}
% It is worth noting that for $\eta_1=\eta_2=1$, the propagation coefficient of the split-step method tends to zero for large values of $\lambda \tau$ meaning that the scheme has a small relaxation time.
% It is also clear that by an appropriate choice of the parameter $\theta$, the scheme is able to reproduce the exact stationary variance of the true solution.
Moreover, by an appropriate choice of the parameter $\theta$, the scheme can reproduce the exact stationary variance of the true solution \cchange{for each given $\lambda\tau$}.
% This can be ensured by the following condition
The following condition can ensure this

\begin{align}\label{eq:theta_equation}
	\frac{2 \lambda\tau}{(1+\theta\lambda\tau)^2 - \displaystyle{\frac{1}{(1+(1-\theta)\lambda\tau)^2}}}=1.
\end{align}
The above equation does not have a closed form solution. % in terms of elementary functions. 
However, for large values of $\lambda\tau$, the fraction in the denominator can be neglected which gives the following approximation

$$\theta = \sqrt{\frac{2}{\lambda\tau}}-\frac{1}{\lambda\tau}, \qquad \lambda=c_1+c_2.$$
For the small values of $\lambda \tau$, one can search for the solution in the form of

\begin{align*}
	\theta &= \theta^{'} + \theta^{''} \lambda \tau + \theta^{'''} ( \lambda \tau )^2 + ...
\end{align*}
Substituting this asymptotic expansion into the original equation and solving for the coefficients, we get

\begin{align*}
	\theta^{'} &= \frac{3-\sqrt{3}}{2}, \quad 
	\theta^{''} = \frac{-9+5\sqrt{3}}{6}, \quad 
	\theta^{'''} = \frac{108 - 187\sqrt{3}}{24}, \quad ...
\end{align*}
As can be seen from Figure \ref{fig:theta_choice}, it is appropriate to choose the following representation for the parameter $\theta$

\begin{align}\label{eq:theta_choice}
	\theta(\lambda \tau) = 
	\begin{cases} 
		\theta^{'} + \theta^{''} \lambda \tau,                & \text{if } \lambda \tau \leq 2.45, \\[1em]
		\sqrt{\frac{2}{\lambda\tau}}-\frac{1}{\lambda\tau},   & \text{if } \lambda \tau > 2.45.
	\end{cases}
\end{align}
\change{The value $\lambda \tau = 2.45$ corresponds to the point where the two branches of $\theta(\lambda \tau)$ coincide.}

\ccchange{
Note that the above formula is obtained by fitting the asymptotical moments of the numerical solution to the exact stationary moments of the true solution.
This approach is suitable only for simple systems with known analytical distributions.
In the next section, we generalize this parameter estimation technique to arbitrary chemical networks.
}

The parameter $\lambda = c_1 + c_2$ in the formula \eqref{eq:theta_choice} is the measure of the relaxation rate of the reversible pair $S_1 \underset{c_2}{\stackrel{c_1}{\rightleftharpoons}} S_2$ to its stationary distribution.
It can also be obtained through a deterministic analysis by studying the eigenvalues of the Jacobian matrix of the corresponding reaction rate equation.
Indeed, the coefficient matrix of the linear differential system

\begin{align*}
	\frac{d}{dt}
	\begin{pmatrix}
		z_1  \\[0.5em]
		z_2 
	\end{pmatrix}
	= 
    %\frac{1}{V}
	\begin{pmatrix}
		-c_1 &  c_2  \\[0.5em]
		 c_1 & -c_2 
	\end{pmatrix}
	\begin{pmatrix}
		z_1  \\[0.5em]
		z_2 
	\end{pmatrix}	
\end{align*}
has two eigenvalues: $\lambda_1=0$ and \cchange{$\lambda_2=-c_1-c_2$}.
In a general case, the system of RREs with an equilibrium point $X^{*}$ can be linearized around this point giving the system of ODEs
\begin{align*}
	\frac{d e}{dt} = 
	\nu \cdot \mathcal{D} a(X^{*}) \cdot e,
\end{align*}
where $e = X - X^{*}$ and $\mathcal{D} a(X^{*})$ is the Jacobian matrix of the vector with reaction propensities.
For the reversible pairs of reactions, the rank of the stoichiometric matrix $\nu$ is equal to one which, according to the rank-nullity theorem, implies that the coefficient matrix of the linearized system has exactly one non-zero eigenvalue.
Since the trace of a matrix is the sum of its eigenvalues, the rate of decay of a reversible pair can be easily computed as
\begin{align}\label{eq:lambda}
	\lambda = - \Tr \big[ \nu \cdot \mathcal{D} a(X^{*}) \big],
\end{align}
where $\Tr$ denotes the trace operator.

\subsection{Parameter estimation for general systems}	
\label{sec:linear_systems}

% \change{
% The formula \eqref{eq:theta_choice} is obtained by fitting the asymptotical moments of the numerical solution to the exact stationary moments of the true solution.
% This approach is suitable only for very simple systems with known analytical distributions.
% In this section, we generalize this parameter estimation technique to arbitrary chemical networks.
% }

Recall that for the case of linear propensity functions in \eqref{eq:linear_prop}, the temporal evolution of the mean and covariance is described by the system of equations % in \eqref{eq:linear_mean_evol}-\eqref{eq:linear_var_evol}

\begin{align*}
	\tag{\ref{eq:linear_mean_evol}}
	\frac{\partial}{\partial t} \E{X(t)} &= \nu C \cdot \E{X(t)} + \nu d,
	\\ \tag{\ref{eq:linear_var_evol}}
	\frac{\partial}{\partial t} \Cov{X(t)} &= \nu C \cdot \Cov{X(t)} + \Cov{X(t)} \cdot (\nu C)^T 
	\\ \nonumber
	&+ \nu \cdot diag\left( C \cdot \E{X(t)} + d \right) \cdot \nu^T.
\end{align*}
\cchange{ We use these equations to estimate the parameters of the split-step scheme by optimizing the error in the first two moments.

The key idea of the proposed approach can be broken down to three steps.
Firstly, to find the reference values of the mean and covariance, we need to solve the system in \eqref{eq:linear_mean_evol}-\eqref{eq:linear_var_evol} either analytically or numerically.
% Firstly, we need to solve the system in \eqref{eq:linear_mean_evol}-\eqref{eq:linear_var_evol} either analytically or numerically.
% Analytical solutions can 
When the analytical solution is not available, we find its numerical approximation by discretizing the above equations as follows % with the classical theta scheme which yields
}

% Denote by $\mu_n$ and $\sigma_n$ the numerical approximations of $\E{X(t_n)}$ and $\Cov{X(t_n)}$ respectively. 
% % The classical theta scheme gives
% The numerical approximation of the above equations with the classical theta scheme gives

\begin{align}
	\label{eq:lyap_mean}
	\mu_{n+1} &= P_1 \mu_n + \tau p_2,
    \\
	\label{eq:lyap_var}
	P_3 \sigma_{n+1}  + \sigma_{n+1} P_3^T 
	&= P_4 \sigma_{n} + \sigma_{n}   P_4^T 
	%\\ \nonumber
	+ \tau \nu \cdot diag( C \mu_{n+1} + d ) \cdot \nu^T,
\end{align}
\cchange{where $\mu_n$ and $\sigma_n$ are the numerical approximations of $\E{X(t_n)}$ and $\Cov{X(t_n)}$ respectively.}
The coefficient matrices and vectors in the above formulas are given by

% \begin{alignat*}{2}
% 	P_1 &= \Big( I - \tau \nu C \Big)^{-1}, \qquad && P_3 = \frac{1}{2} I - \tau \nu C, 
% 	\\
% 	p_2 &= \Big( I - \tau \nu C \Big)^{-1} \nu d , \qquad && P_4 = \frac{1}{2} I .
% \end{alignat*}
\begin{alignat*}{2}
	P_1 &= \Big( I - \tilde{\theta}\tau \nu C \Big)^{-1} \Big( I + (1-\tilde{\theta})\tau \nu C \Big), \qquad && P_3 = \frac{1}{2} I - \tilde{\theta}\tau \nu C, 
	\\
	p_2 &= \Big( I - \tilde{\theta}\tau \nu C \Big)^{-1} \nu d , \qquad && P_4 = \frac{1}{2} I + (1-\tilde{\theta})\tau \nu C .
\end{alignat*}
\cchange{Note that the parameter $\tilde{\theta}$ above is in general not related to the parameter $\theta$ in the considered split-step scheme \eqref{eq:mod_split_step}.
% In fact, 
Any other approximation scheme for the solution of \eqref{eq:linear_mean_evol}-\eqref{eq:linear_var_evol} can also be used instead of \eqref{eq:lyap_mean}-\eqref{eq:lyap_var}.
Moreover, in all numerical simulations below, we used the fully implicit scheme with $\tilde{\theta} = 1$. }

\cchange{Secondly, we explicitly derive the equations for the mean and covariance of the solution generated with the split-step scheme \eqref{eq:mod_split_step}.
In the case of linear propensity functions, it reads as }

% Similarly, the split-step scheme \eqref{eq:mod_split_step} with the linear propensities in \eqref{eq:linear_prop} can be written as

\begin{align}\label{eq:lin_split_step}
	\nonumber
	\hat{Y}_{n}	&= R_1 Y_n + \tau r_2,
	\\
	\tilde{Y}_{n} &= \hat{Y}_{n} + \nu \mathcal{P} \Big( \big(C \hat{Y}_n+d\big) \tau \Big) - \tau \Big( \nu C \hat{Y}_n + \nu d \Big),
	\\ \nonumber
	Y_{n+1}	&= R_3 \tilde{Y}_{n} + \tau r_4,
\end{align}	
where

\begin{align*}
	R_1 &= \Big( I - \tau \nu \cdot diag(\eta_1) \cdot diag(1-\theta) \cdot C \Big)^{-1} \Big( I + \tau \nu \cdot diag(1-\eta_1) \cdot diag(1-\theta) \cdot C \Big),
	\\
	r_2 &= \Big( I - \tau \nu \cdot diag(\eta_1) \cdot diag(1-\theta) \cdot C \Big)^{-1} \nu \cdot diag(1-\theta) \cdot d,
	\\
	R_3 &= \Big( I - \tau \nu \cdot diag(\eta_2) \cdot diag(\theta) \cdot C \Big)^{-1} \Big( I + \tau \nu \cdot diag(1-\eta_2) \cdot diag(\theta) \cdot C \Big),
	\\
	r_4 &= \Big( I - \tau \nu \cdot diag(\eta_2) \cdot diag(\theta) \cdot C \Big)^{-1} \nu \cdot diag(\theta) \cdot d
\end{align*}
and $\theta=[\theta_1,...,\theta_R]^T$, $\eta_1=[\eta_{11},...,\eta_{1R}]^T$, $\eta_2=[\eta_{21},...,\eta_{2R}]^T$.

By applying Lemma \ref{def:law_total_mean}, it is easy to show that the expectation of $Y_{n+1}$ is given by

\begin{align}
	\nonumber
	\E{\hat{Y}_{n}} &= R_1 \E{Y_n} + \tau r_2,
	\\[0.5em] \label{eq:split_step_mean}
	\E{{Y}_{n+1}} &= R_3 \E{\hat{Y}_{n}} + \tau r_4.
\end{align}
To find the covariance of $Y_{n+1}$, recall that the components of the random vector $\mathcal{P} \Big( \big(C \hat{Y}_n+d\big) \tau \Big)$ are \cchange{ independent Poisson variables conditioned on $Y_n$ and thus $\Cov{\mathcal{P} \Big( \big(C \hat{Y}_n+d\big) \tau \Big)|Y_n}$ is a diagonal matrix}.  
This gives

\begin{align}
	\nonumber
	\Cov{\hat{Y}_{n}} &= R_1 \Cov{Y_n} R_1^T ,
	\\[0.5em] \label{eq:split_step_var}
	\Cov{\tilde{Y}_{n}} &= \Cov{\hat{Y}_{n}} + \tau \nu \cdot diag \left( C \E{\hat{Y}_n} + d \right) \cdot \nu^T,
	\\[0.5em] \nonumber
	\Cov{Y_{n+1}} &= R_3 \Cov{\tilde{Y}_{n}} R_3^T.
\end{align}

\cchange{Finally, by comparing the equations for the means in \eqref{eq:linear_mean_evol}, \eqref{eq:split_step_mean} and covariances in \eqref{eq:linear_var_evol}, \eqref{eq:split_step_var}, we obtain the optimal parameter values $\theta^n$, $\eta_{1}^n$ and $\eta_{2}^n$ at each time step as a solution to the following minimization problem }

% The matrices $R_1$, $R_3$ and the vectors $r_2$, $r_4$ in \eqref{eq:lin_split_step} are the functions of parameters $\theta$, $\eta_{1}$ and $\eta_{2}$.
% By comparing the equations for the means in \eqref{eq:linear_mean_evol}, \eqref{eq:split_step_mean} and covariances in \eqref{eq:linear_var_evol}, \eqref{eq:split_step_var}, the optimal parameter values at each time step can be obtained as the solution to the following minimization problem

\begin{align}\label{eq:argmin}
	\{ \theta^n, \eta_1^n, \eta_2^n \} = \argmin_{ \theta^n,\eta_{1}^n,\eta_{2}^n \in [a,b] } \Big( \norm{\E{Y_{n}}-\mu_{n}}{}^2 + \norm{\Cov{Y_{n}}-\sigma_{n}}{f}^2 \Big) ,
\end{align}
where $[a,b]$ is a bounded interval, $\norm{\cdot}{f}$ is the Frobenius matrix norm and $\norm{\cdot}{}$ is the usual Euclidean vector norm.
Algorithm \ref{alg:scheme} contains a detailed description of the proposed parameter estimation approach.

\renewcommand{\algorithmicrequire}{\textbf{Input:}}
\renewcommand{\algorithmicensure}{\textbf{Output:}}
\begin{algorithm}[t]
\caption{Parameter estimation for the split-step tau-leaping algorithm}
\label{alg:scheme}
\begin{algorithmic}[1]
	\REQUIRE{Initial condition $X_0$, time step $\tau$, final time $T$, $\alpha_1 < 1$, $\alpha_2 < 1$}
	\ENSURE{Optimal parameter values}
%	\STATE{	Choose $\theta_0 \in [0,1]$ in \eqref{eq:lyap_mean}-\eqref{eq:lyap_var} }
%	\STATE{	Choose $\alpha_1 < 1$ and $\alpha_2 < 1$ }
	\FOR{ each reversible reaction pair }
		\STATE{Estimate the relaxation rate $\lambda_r$ using formula \eqref{eq:lambda}}
		\STATE{Find $\theta_r^0$ as a solution of the equation \eqref{eq:theta_equation} 	or using formula \eqref{eq:theta_choice} }
%		\begin{align*}
%			(1+(1-\theta_r^0)\lambda_r\tau)^2 \Big( 2 \lambda_r\tau - (1+\theta_r^0\lambda_r\tau)^2 \Big) = 1
%%			2 \lambda\tau = (1+\theta_r^0\lambda_r\tau)^2 - \displaystyle{\frac{1}{(1+(1-\theta_r^0)\lambda_r\tau)^2}}
%%			\theta_r^0(\lambda_r \tau) = 
%%			\begin{cases} 
%%				\Big( 9-3\sqrt{3} + \big(-9+5\sqrt{3}\big) \lambda_r \tau \Big) / 6,                     & \text{if } \lambda_r \tau \leq 2, \\[1em]
%%				\sqrt{\frac{2}{\lambda_r\tau}}-\frac{1}{\lambda_r\tau},   & \text{if } \lambda_r \tau > 2
%%			\end{cases}
%		\end{align*}
		\STATE{Set $\eta_{1r}^0=\eta_{2r}^0=1$}
	\ENDFOR
	%\ENSURE{Find optimal parameters}
	\STATE{$n \leftarrow 1$}
	\STATE{$t_1=t_{0}+\tau$}
	\WHILE{$t_n \leq T$}
		\STATE{Find $\mu_n$ and $\sigma_n$ from \eqref{eq:lyap_mean} and \eqref{eq:lyap_var}}
		\STATE{Set $C=\mathcal{D} a\big(\mu_n\big)$ and $d=a\big( \mu_n \big) - C \cdot \mu_n$}
		\STATE{Set $\theta^{n-1}, \eta_1^{n-1}, \eta_2^{n-1}$ as initial guess for $\theta^n, \eta_1^n, \eta_2^n$}
		\STATE{Find $\theta^n, \eta_1^n, \eta_2^n$ from}
		\begin{align*}
			\{ \theta^n, \eta_1^n, \eta_2^n \} = \argmin \Big( \norm{\E{Y_{n}}-\mu_{n}}{}^2 + \norm{\Cov{Y_{n}}-\sigma_{n}}{f}^2 \Big)
		\end{align*}
		\STATE{with $\E{Y_n}$ and $\Cov{Y_n}$ given by the formulas \eqref{eq:split_step_mean} and \eqref{eq:split_step_var} respectively}
		\IF{ $\frac{\norm{\E{Y_{n}}-\mu_{n}}{}}{\norm{\mu_{n}}{}} \geq \alpha_1$ \OR $\frac{\norm{\Cov{Y_{n}}-\sigma_{n}}{f}}{\norm{\sigma_{n}}{f}} \geq \alpha_2$ }
			\STATE{Reduce $\tau$}
		\ELSE
			\STATE{$n \leftarrow n+1$}
		\ENDIF
		\STATE{$t_n=t_{n-1}+\tau$}
	\ENDWHILE
\end{algorithmic}
\end{algorithm}

For nonlinear systems, one can use the following linearization %can be used

\begin{align}\label{eq:linearized_prop}
	a(X) 
	&\approx a\big( \E{X} \big) + \mathcal{D} a\big(\E{X}\big) \cdot \big( X - \E{X} \big).
\end{align}
By setting 

$$C=\mathcal{D} a\big(\E{X}\big) \qquad \text{and} \qquad d=a\big( \E{X} \big) - C \cdot \E{X},$$
this approximation allows applying Algorithm \ref{alg:scheme} to arbitrary chemical networks.

\textbf{Remark.} Special care must be taken when addressing the existence and uniqueness of the solution to the optimization problem in \eqref{eq:argmin}.
The uniqueness part can be relaxed because any sufficiently small local minimum of the target functional can be taken as a solution.
Also, the initial guess with $\eta_{1}=\eta_{2}=1$ and $\theta$ as in \eqref{eq:theta_choice} (by considering reversible pairs individually in isolation from other reactions) often lies in the basin of attraction of the global minimum.
The existence of the solution is also guaranteed because the problem is constrained to the bounded interval.

\section{Numerical examples}
\label{sec:num_examples}

% The goal of this section is to give a numerical justification of the proposed splitting approach.

% We present the numerical results for both linear and nonlinear problems with 

\change{
%The numerical results in this section show that the proposed split-step method is well suited for the integration of fast stochastic reactions.
In this section, we present numerical results for both linear and nonlinear problems.
The provided examples aim to show that the proposed integration technique is capable to accurately recover stationary distributions of reaction networks without explicitly resolving their temporal scales.
% well-suited for the integration of the fast and stable fast reactions 
% which show that the proposed split-step tau-leap method is well suited for the integration of fast stochastic reactions. % and compare its performance to the classical implicit theta schemes.
%This section contains numerical results for both linear and nonlinear problems.
% In all examples below, the estimation of parameters was performed with the single run of the standard unconstrained gradient search with the initial conditions as in Algorithm \ref{alg:scheme}.
In all examples below, we estimated the parameters with the single run of the standard unconstrained gradient search with the initial conditions as in Algorithm \ref{alg:scheme}.
All simulations were performed using $10^6$ Monte Carlo samples.
}

\subsection{Example 1. (Linear system)}

Consider the monomolecular reaction network with three pairs of reversible reactions
\begin{align}\label{ex:1}
S_1~\underset{c_2}{\stackrel{c_1}{\rightleftharpoons}}~S_2~\underset{c_4}{\stackrel{c_3}{\rightleftharpoons}}~S_3~\underset{c_6}{\stackrel{c_5}{\rightleftharpoons}}~S_4.
\end{align}
This network is described by the following stoichiometric matrix

\begin{align*}
	\nu = 
	\begin{pmatrix}
		-1 &  1 &  0 &  0 &  0 &  0 \\[0.5em]
		 1 & -1 & -1 &  1 &  0 &  0 \\[0.5em]
		 0 &  0 &  1 & -1 & -1 &  1 \\[0.5em]
		 0 &  0 &  0 &  0 &  1 & -1
	\end{pmatrix}
\end{align*}
and the corresponding propensity functions 

\begin{alignat*}{2}
	&\alpha_1(X(t)) = c_1 X_1(t), \quad & \alpha_2(X(t)) = c_2 X_2(t), & \\
	&\alpha_3(X(t)) = c_3 X_2(t), \quad & \alpha_4(X(t)) = c_4 X_3(t), & \\
	&\alpha_5(X(t)) = c_5 X_3(t), \quad & \alpha_6(X(t)) = c_6 X_4(t). & 
\end{alignat*}

Monomolecular reactions provide an example of chemical systems with known analytical solutions \cite{Jahnke2007}.
For the initial condition $X(0)=[x_T,0,...,0]^T$, the exact distribution of $X(t)$ is multinomial and has the particularly simple form 

\begin{align*}
	\mathbb{P}_{t,x} &= \mathcal{M} \big( x,x_T,p(t) \big),
	\qquad
	p(t) = \exp(t \nu C) \cdot [1,0,...,0]^T. %\varepsilon_1,
\end{align*}
The exact mean and covariance are thus given by

\begin{align*}
	\E{X(t)} &= x_T \cdot p(t),
	\\
	\Cov{X_i(t),X_j(t)} &= 
	\begin{cases} 
		x_T \cdot p_i(t) (1-p_i(t)) & \text{if } i = j, \\
		x_T \cdot p_i(t) p_j(t) & \text{if } i \neq j.
	\end{cases}
\end{align*}

\begin{figure}[t!]
	\centering
	    \begin{subfigure}[t]{1.0\textwidth}
                \centering
                \includegraphics[width=0.45\textwidth]{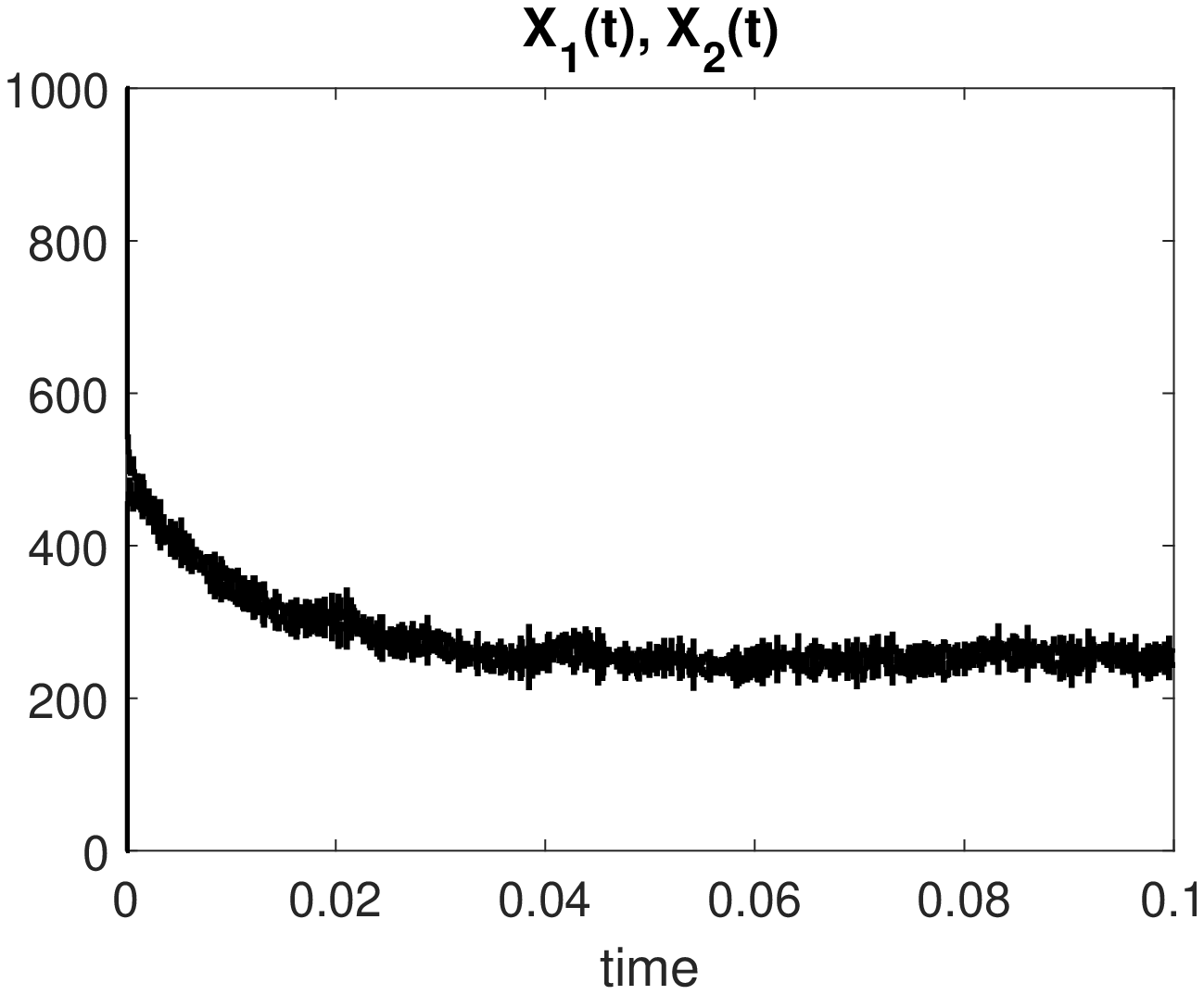}
                \qquad
                \includegraphics[width=0.45\textwidth]{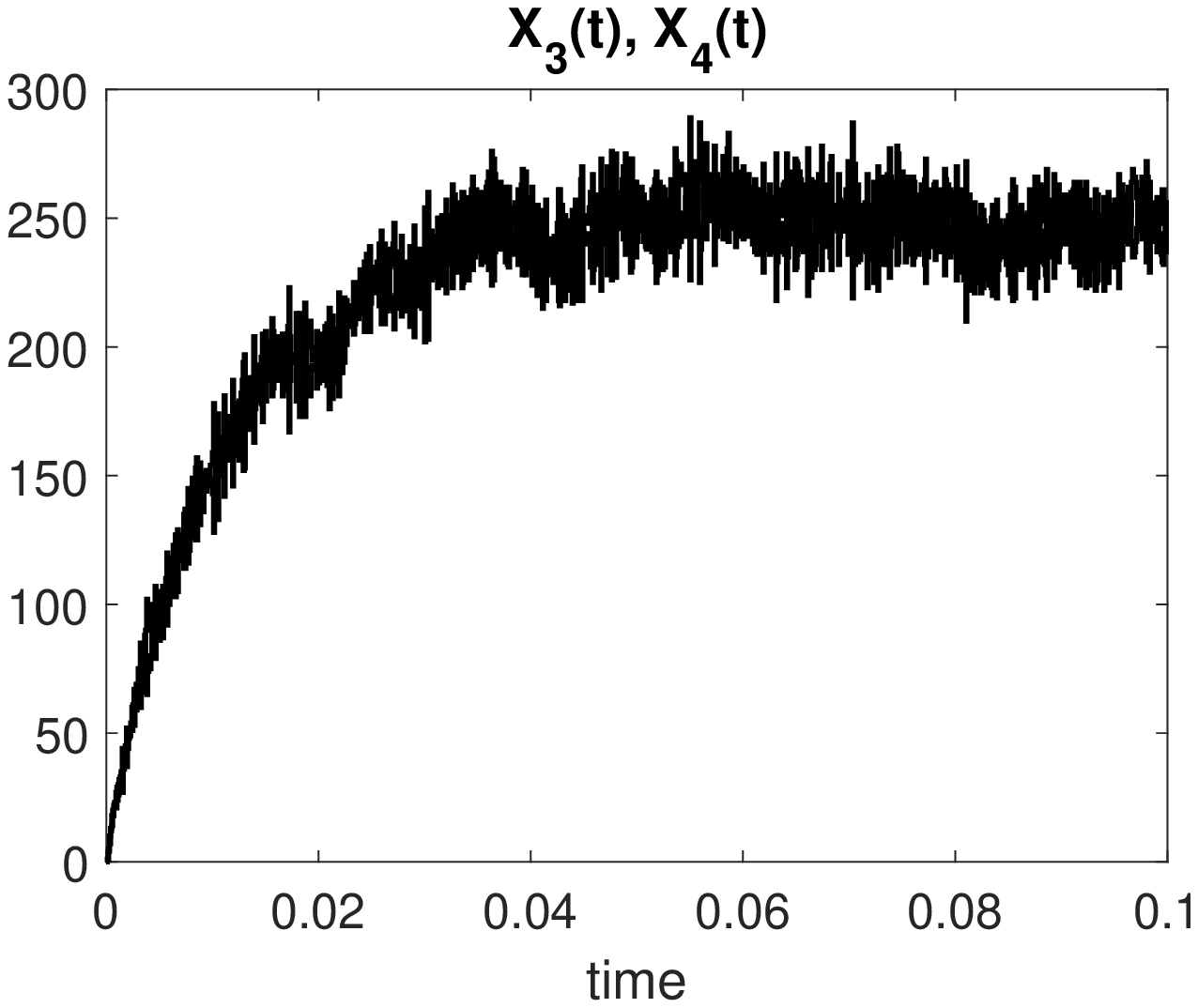}
        \end{subfigure}
	\caption{\cchange{A realization of the state paths of the system in \eqref{ex:1} with the initial condition $X(0)=[1000,0,0,0]^T$ and reaction rates $c~=~[10^4,10^4,10^2,10^2,10^5,10^5]$.}}
	\label{fig:ex_1_state_paths}
\end{figure}

\begin{figure}[t!]
	\centering
	   % \begin{subfigure}[t]{1.0\textwidth}
    %             \centering
    %             \includegraphics[width=0.45\textwidth]{ex_1_xT_1000_X12.eps}
    %             \qquad
    %             \includegraphics[width=0.45\textwidth]{ex_1_xT_1000_X34.eps}
    %     \end{subfigure}
    %     \\
        \begin{subfigure}[t]{1.0\textwidth}
                \centering
                \includegraphics[width=0.45\textwidth]{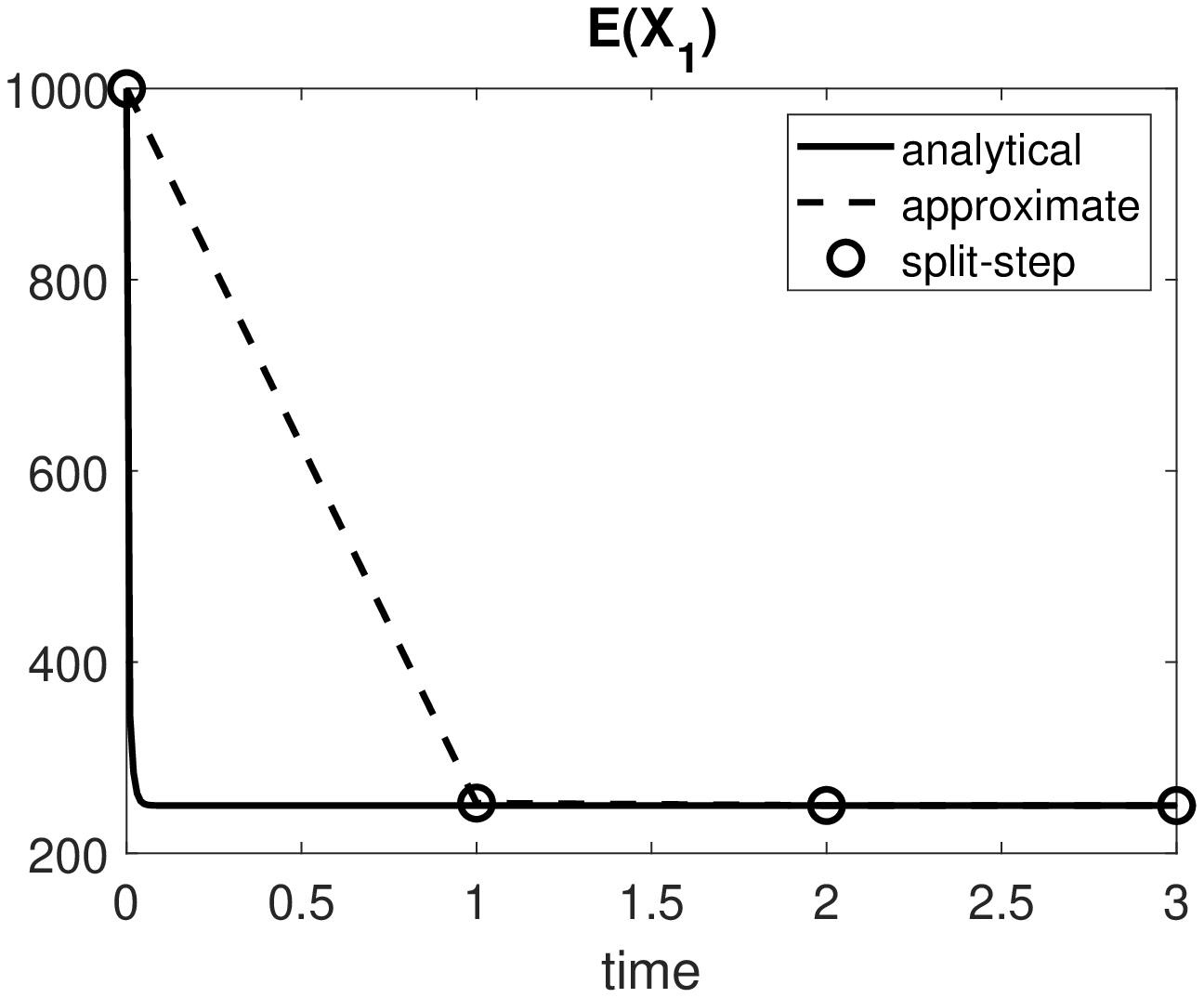}
                \qquad
                \includegraphics[width=0.45\textwidth]{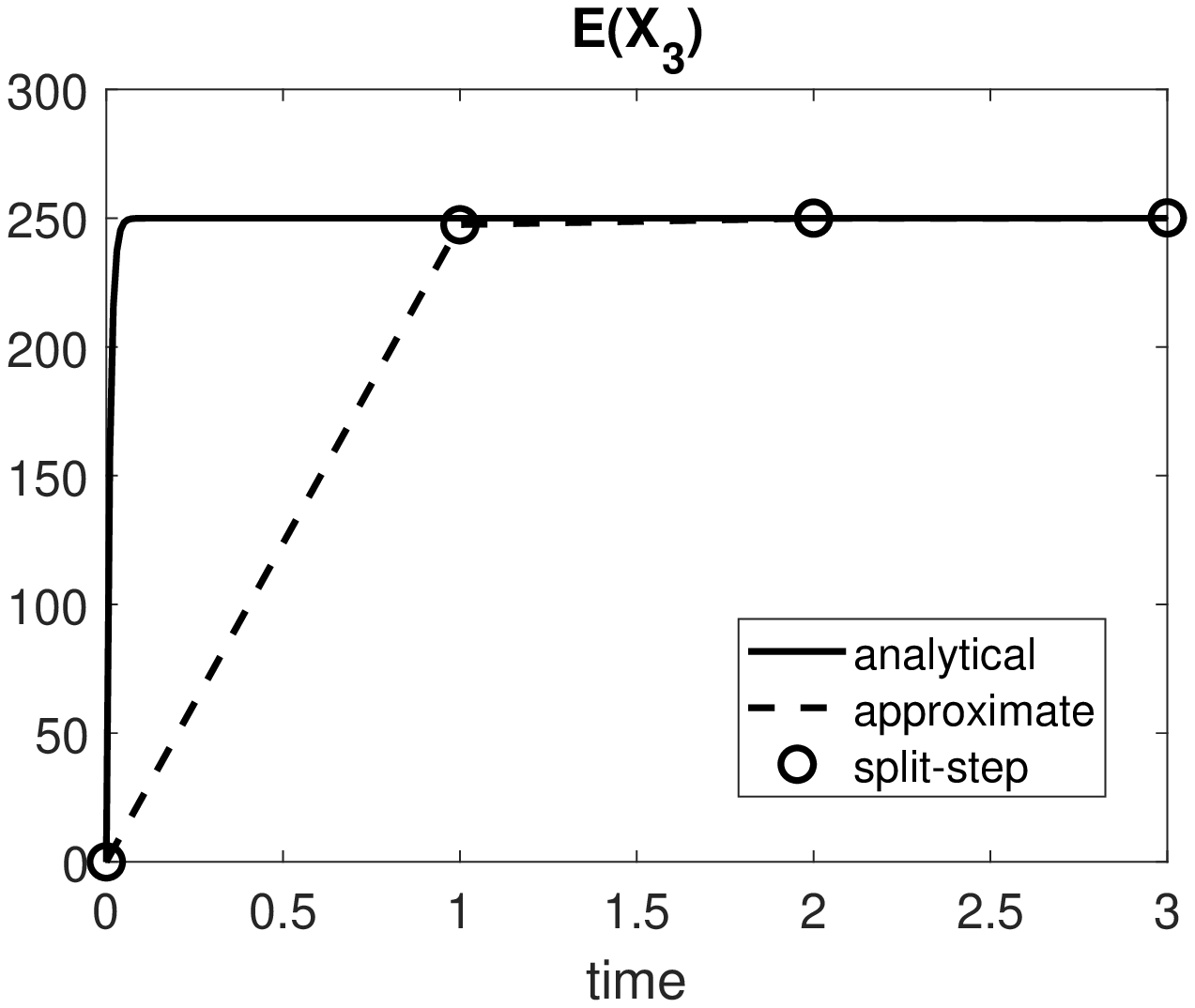}
        \end{subfigure}
        \\
        \begin{subfigure}[t]{1.0\textwidth}
                \centering
                \includegraphics[width=0.45\textwidth]{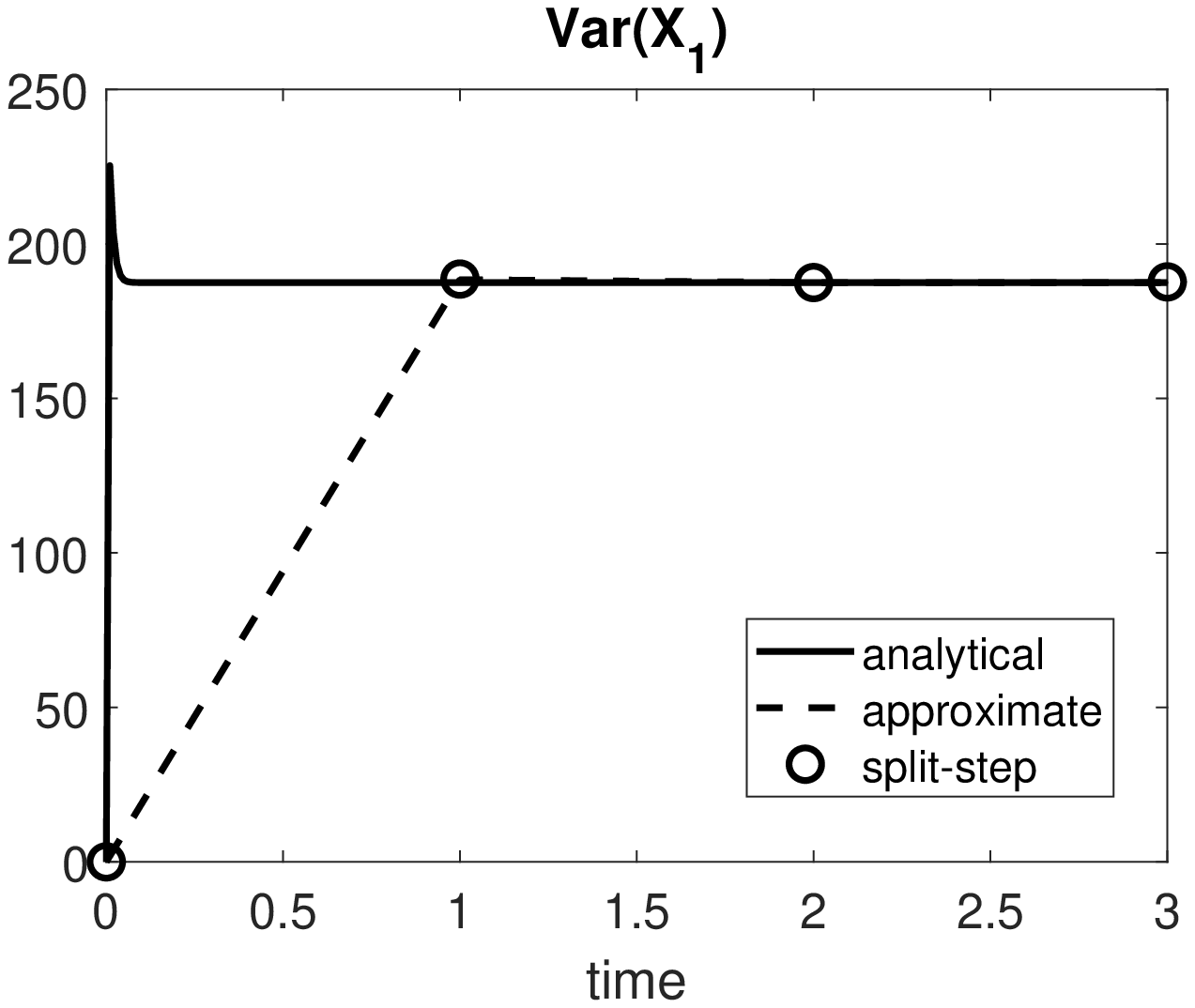}
                \qquad
                \includegraphics[width=0.45\textwidth]{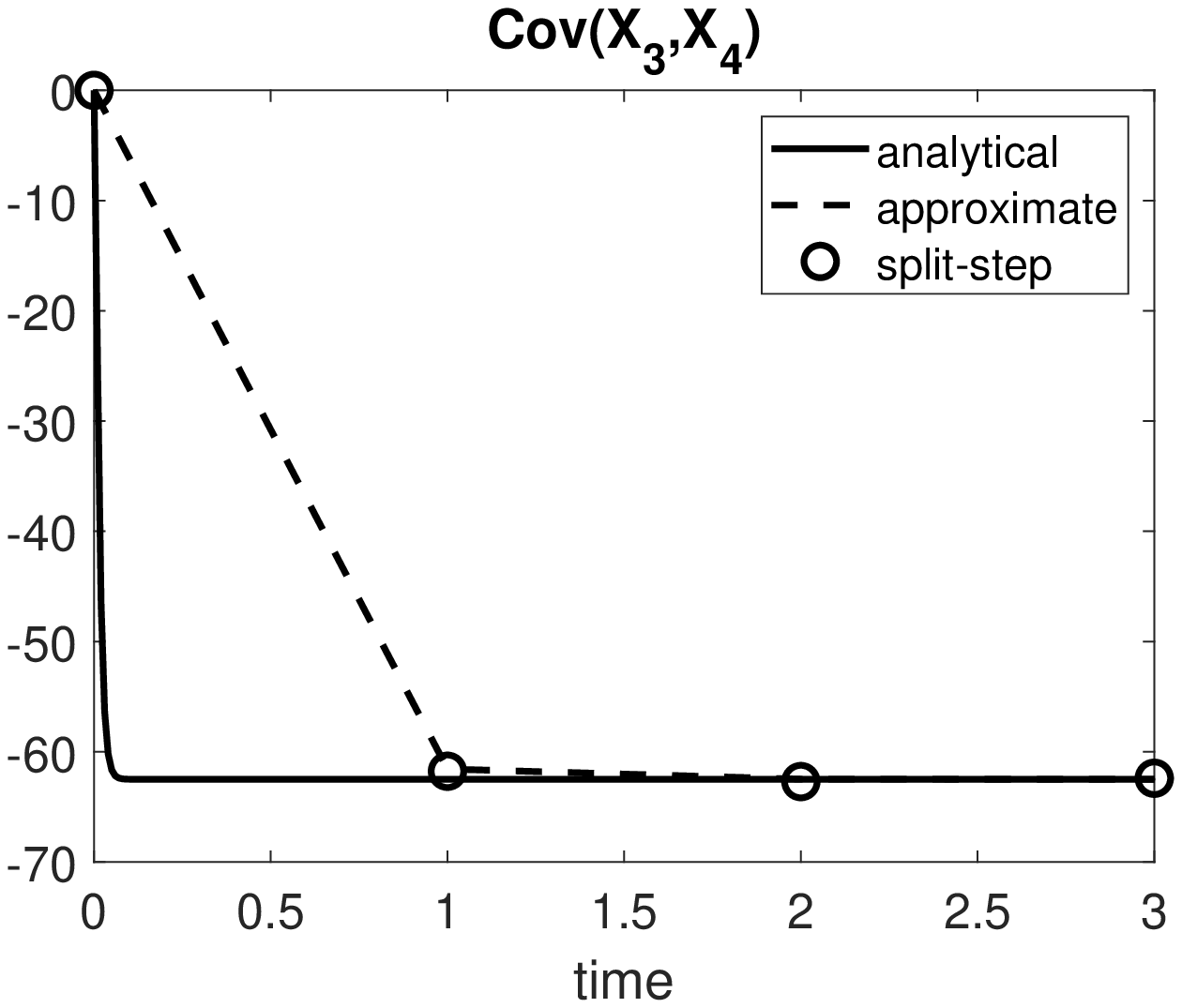}
        \end{subfigure}
	\caption{\cchange{Evolution of the mean and covariance of the system in \eqref{ex:1} with the initial condition $X(0)=[1000,0,0,0]^T$ and reaction rates $c~=~[10^4,10^4,10^2,10^2,10^5,10^5]$.}}
	\label{fig:ex_1_mean_cov_b}
\end{figure}

% \begin{figure}[t!]
% 	\centering
%         \begin{subfigure}[t]{1.0\textwidth}
%                 \centering
%                 \includegraphics[width=0.45\textwidth]{ex_1_xT_1000_c_1e-1_mean1}
%                 \qquad
%                 \includegraphics[width=0.45\textwidth]{ex_1_xT_1000_c_1e-1_mean3}
%         \end{subfigure}
%         \\
%         \begin{subfigure}[t]{1.0\textwidth}
%                 \centering
%                 \includegraphics[width=0.45\textwidth]{ex_1_xT_1000_c_1e-1_var1}
%                 \qquad
%                 \includegraphics[width=0.45\textwidth]{ex_1_xT_1000_c_1e-1_cov34}
%         \end{subfigure}
% 	\caption{Evolution of the mean and covariance of the system in \eqref{ex:1} with the initial condition $X(0)=[1000,0,0,0]^T$ and reaction rates $c~=~10^{-1} \cdot [10,1,10,10,1,5]$.}
% 	\label{fig:ex_1_mean_cov_b}
% \end{figure}

\cchange{Figure \ref{fig:ex_1_state_paths} shows a single realization of the state paths for the stiff system in \eqref{ex:1}.}
The temporal evolution of some of the moments is depicted in Figure \ref{fig:ex_1_mean_cov_b}.
It also shows the approximate moments calculated from the systems in \eqref{eq:lyap_mean}-\eqref{eq:lyap_var} and \eqref{eq:split_step_mean}-\eqref{eq:split_step_var}. 
It is seen that the parameter estimation approach in \eqref{eq:argmin} works well and both approximate solutions tend to the stationary state of the true analytical solution as desired.

\begin{figure}[!t]
	\centering
        \begin{subfigure}[t]{1.0\textwidth}
                \centering
                \includegraphics[width=0.45\textwidth]{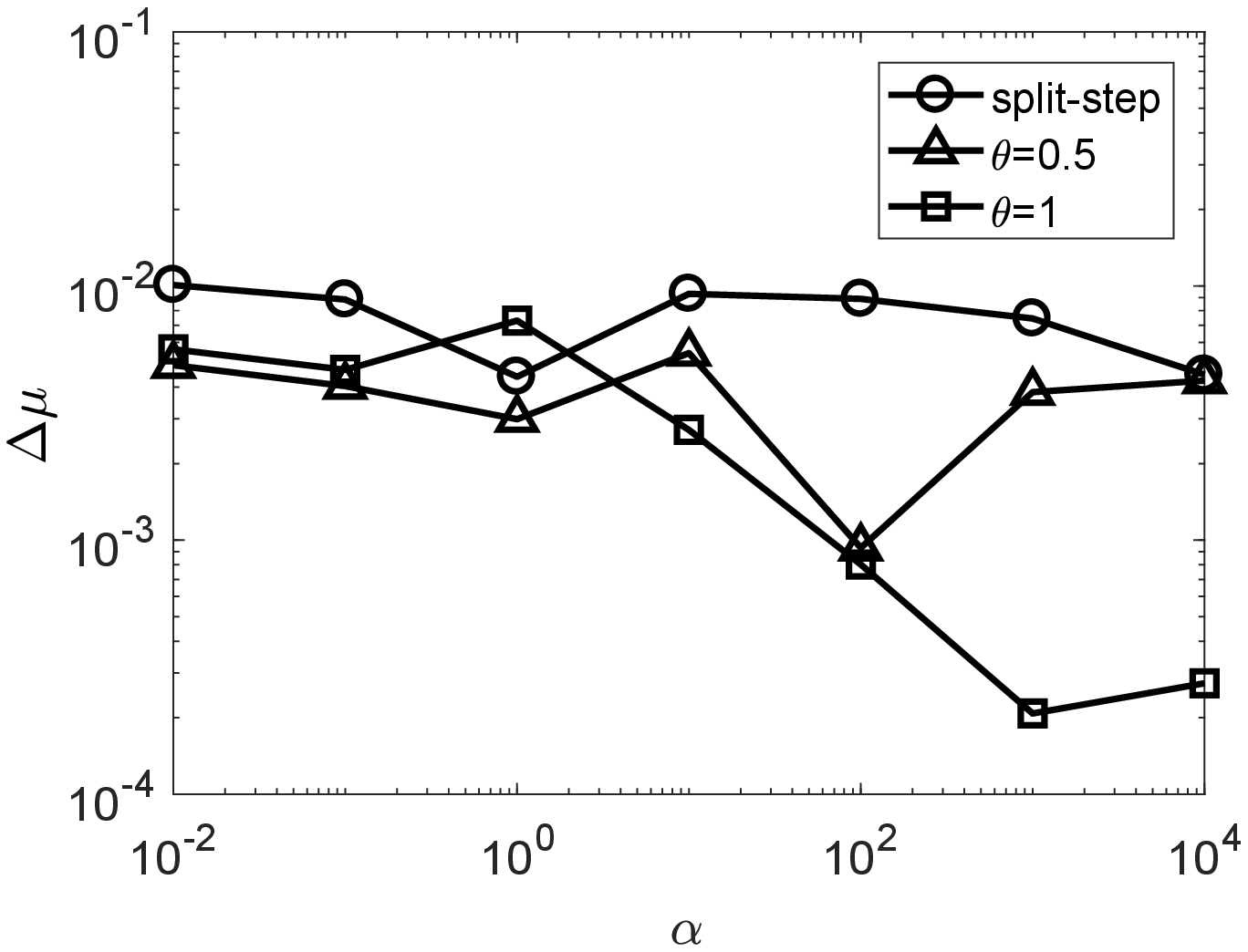}
                \qquad
                \includegraphics[width=0.45\textwidth]{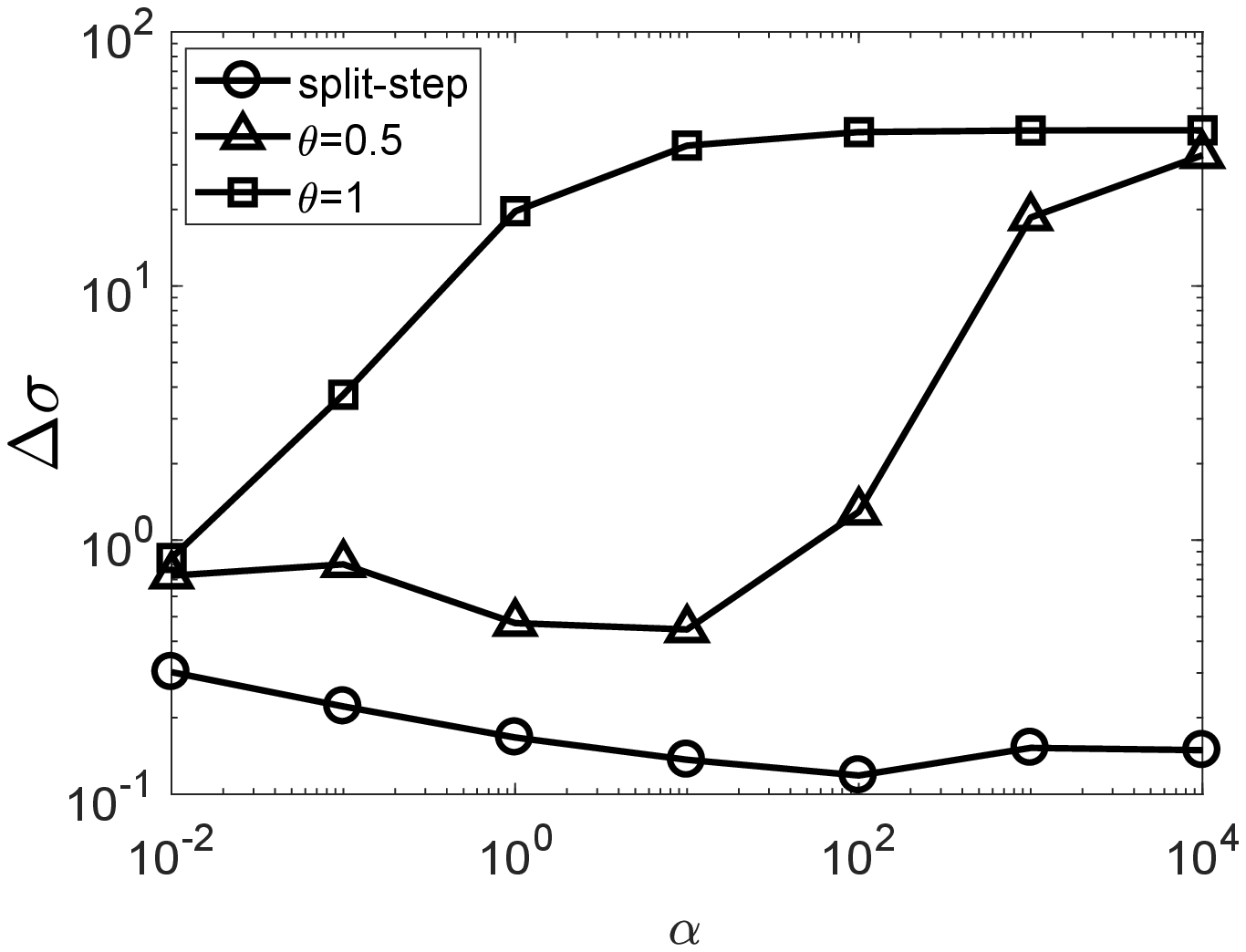}
        \end{subfigure}
	\caption{Errors of the mean (left) and the covariance (right) of the simulated approximate solutions at $T=100$ for the reaction rates $c~=\alpha \cdot [1,1,1,1,1,1]$.}
	\label{fig:ex_1_error_cost}
\end{figure}

To test the accuracy of the integrators, we fix the time step $\tau=1$ and the final time $T=100$ and consider the reaction rates of the form
\begin{align*}
	c = \alpha \cdot [1,1,1,1,1,1]^T.
\end{align*}
% Different values of $\alpha$ correspond to different temporal scales and allow to compare the performance of the proposed and classical methods for different levels of stiffness.
\cchange{Different values of $\alpha$ allow to compare the performance of the proposed and classical methods at different temporal scales.}
Figure \ref{fig:ex_1_error_cost} illustrates the results of this comparison for the split-step algorithm in \eqref{eq:spli_step} and the theta scheme in~\eqref{eq:theta}.
It is clear that both methods give accurate results for the stationary mean.
However, the ability of the theta scheme to preserve the stationary variance is deteriorating \cchange{at faster scales} while the split-step method demonstrates the same level of accuracy uniformly in $\alpha$.
Also, note that the proposed algorithm has two implicit steps and is more expensive than the theta scheme using the same temporal discretization.
However, the results in Figure \ref{fig:ex_1_error_cost} indicate that the split-step method is much more accurate and is more appropriate for simulating fast reactions.

\begin{figure}[!t]
	\centering
        \begin{subfigure}[t]{1.0\textwidth}
                \centering
                \includegraphics[width=0.45\textwidth]{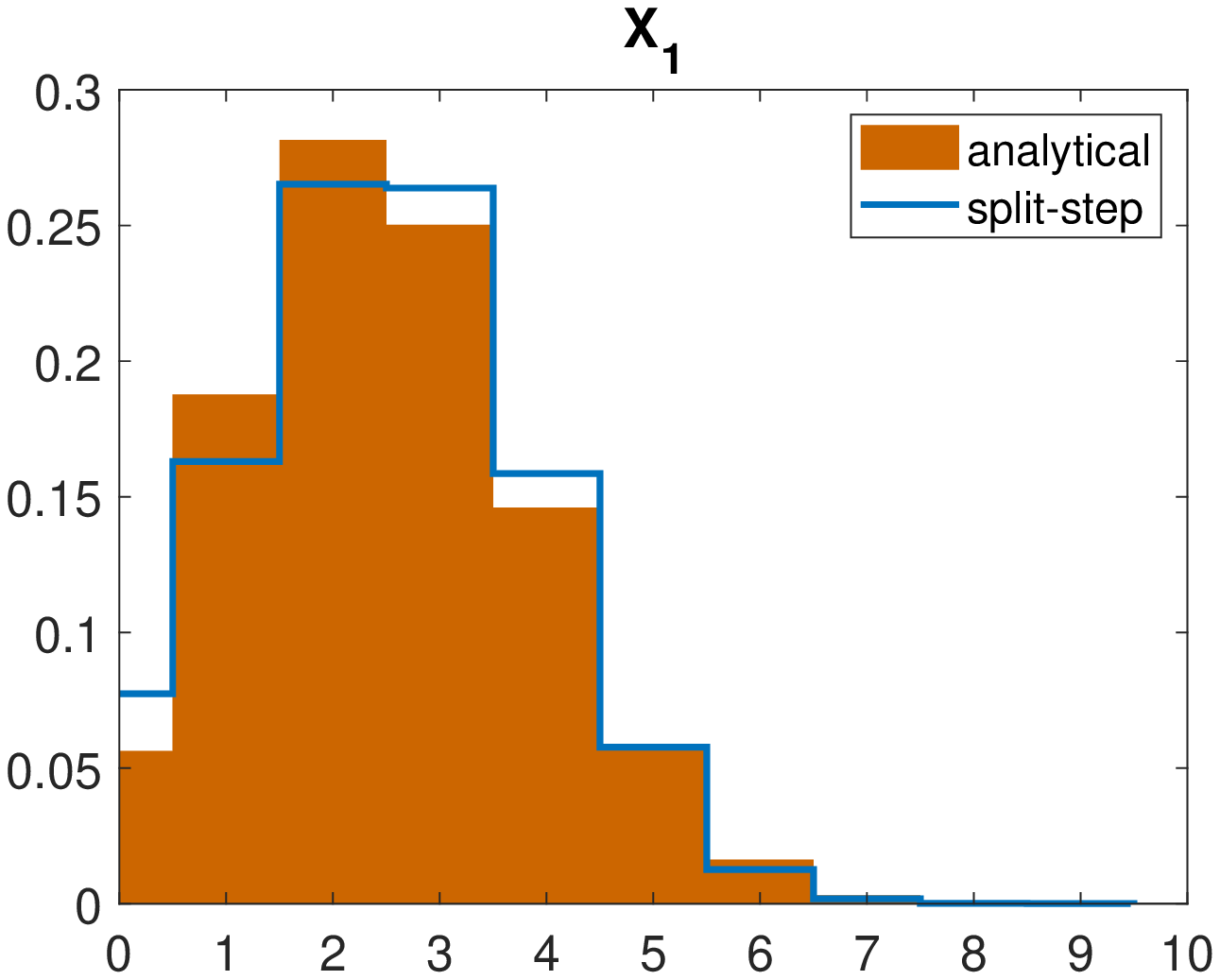}
                \qquad
                \includegraphics[width=0.45\textwidth]{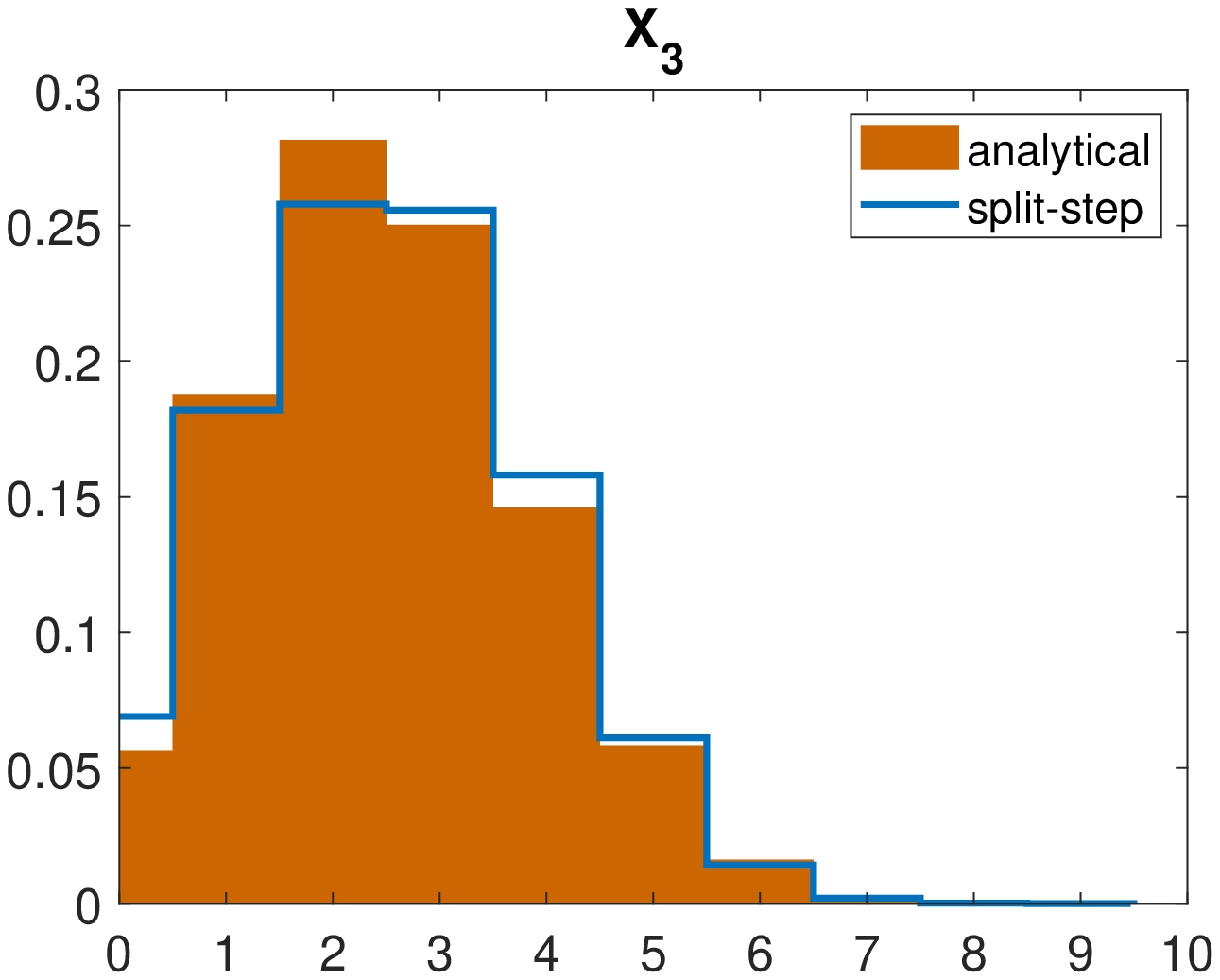}
                \caption{$x_T=10$}
                \label{fig:ex_1_distr_xT=10}
        \end{subfigure}
        \\
        \begin{subfigure}[t]{1.0\textwidth}
                \centering
                \includegraphics[width=0.45\textwidth]{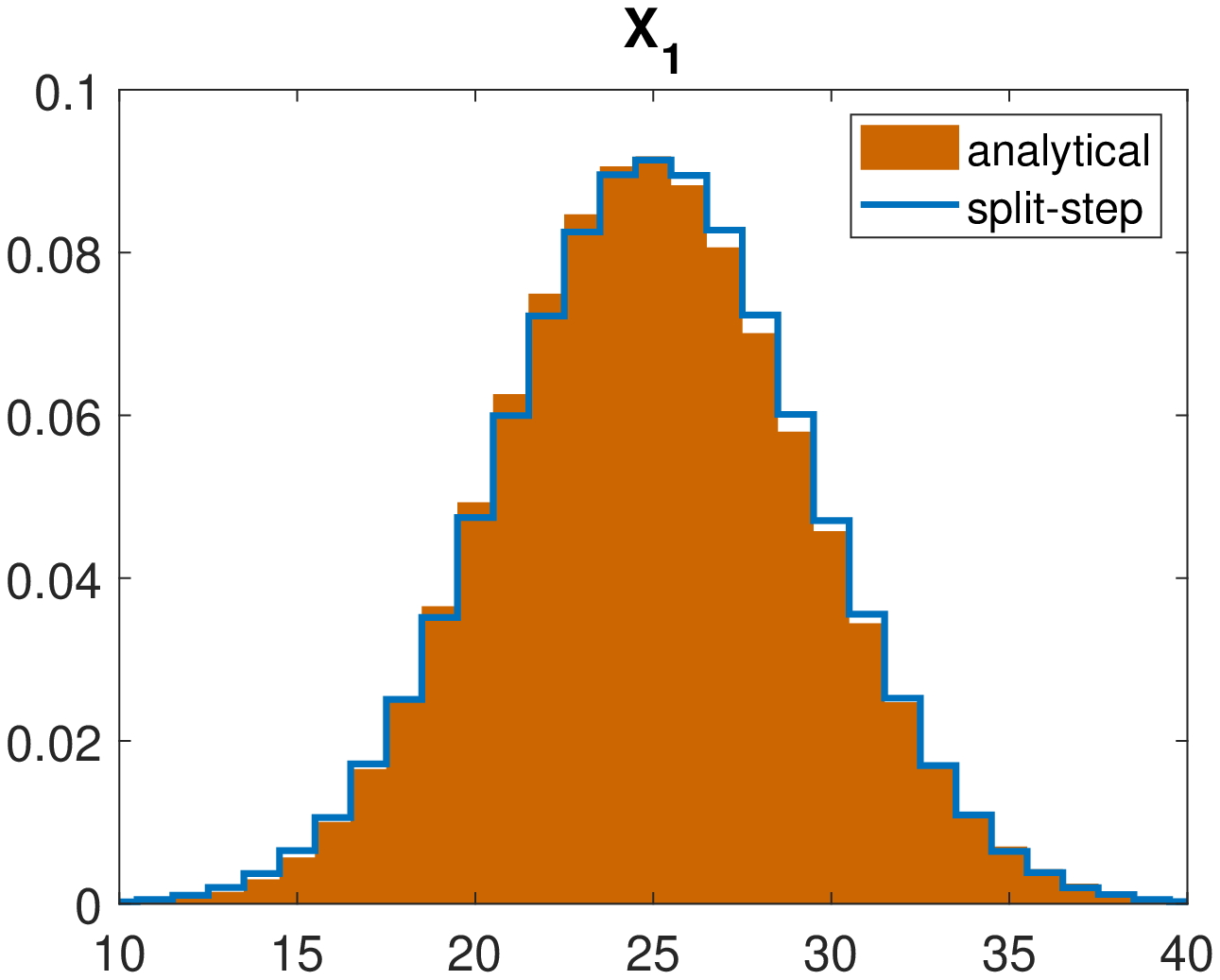}
                \qquad
                \includegraphics[width=0.45\textwidth]{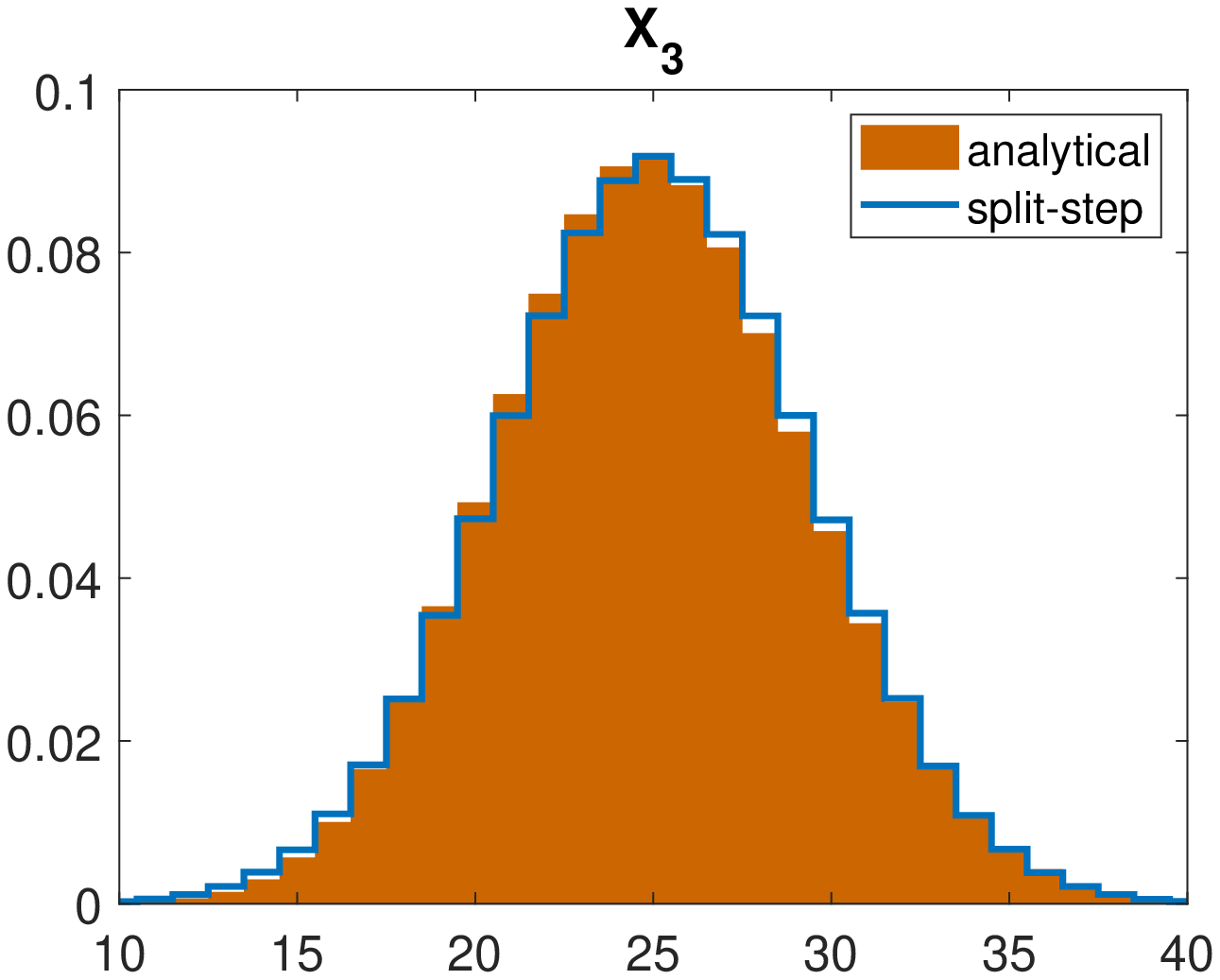}
                \caption{$x_T=100$}
                \label{fig:ex_1_distr_xT=100}
        \end{subfigure}
        \\
        \begin{subfigure}[t]{1.0\textwidth}
                \centering
                \includegraphics[width=0.45\textwidth]{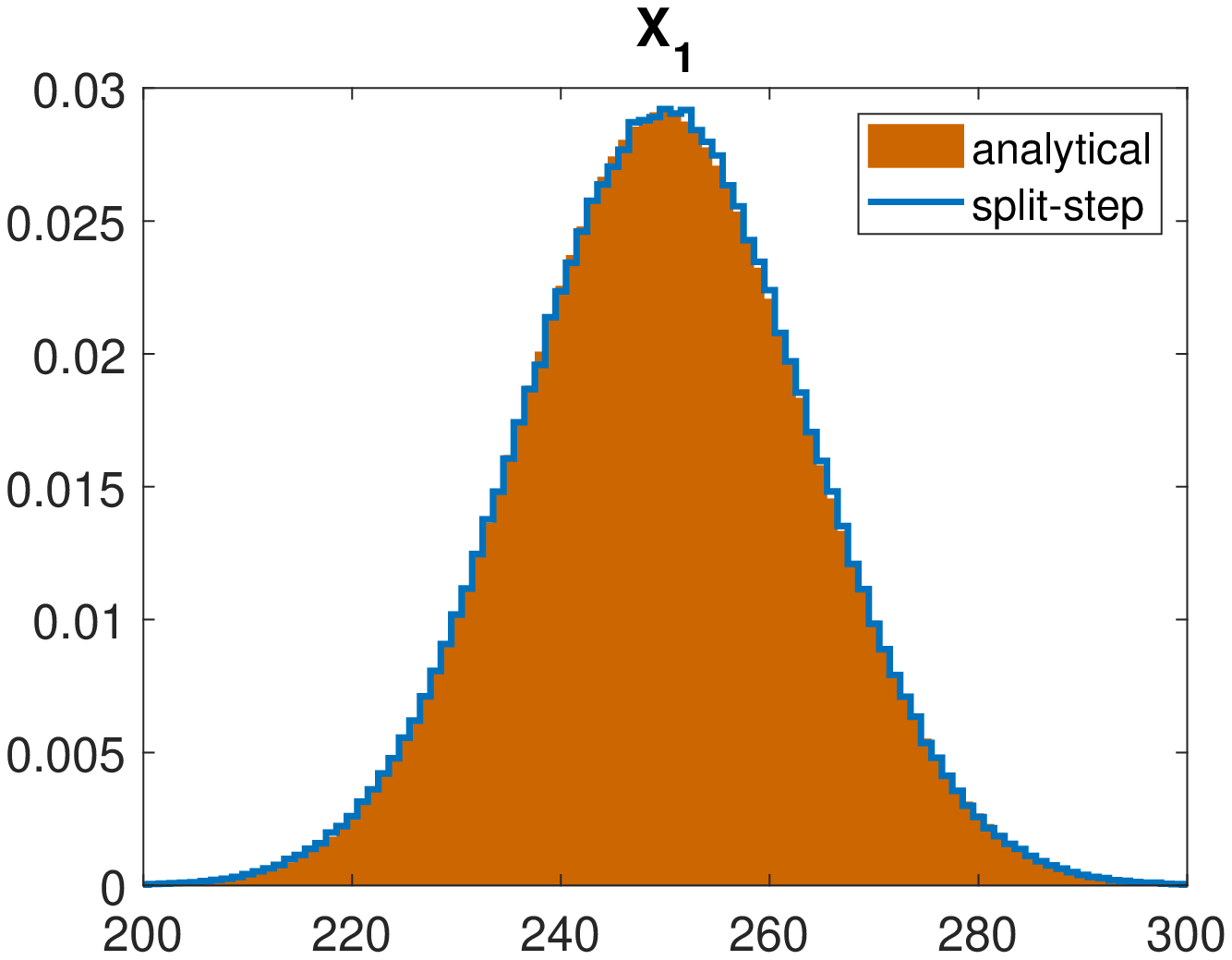}
                \qquad
                \includegraphics[width=0.45\textwidth]{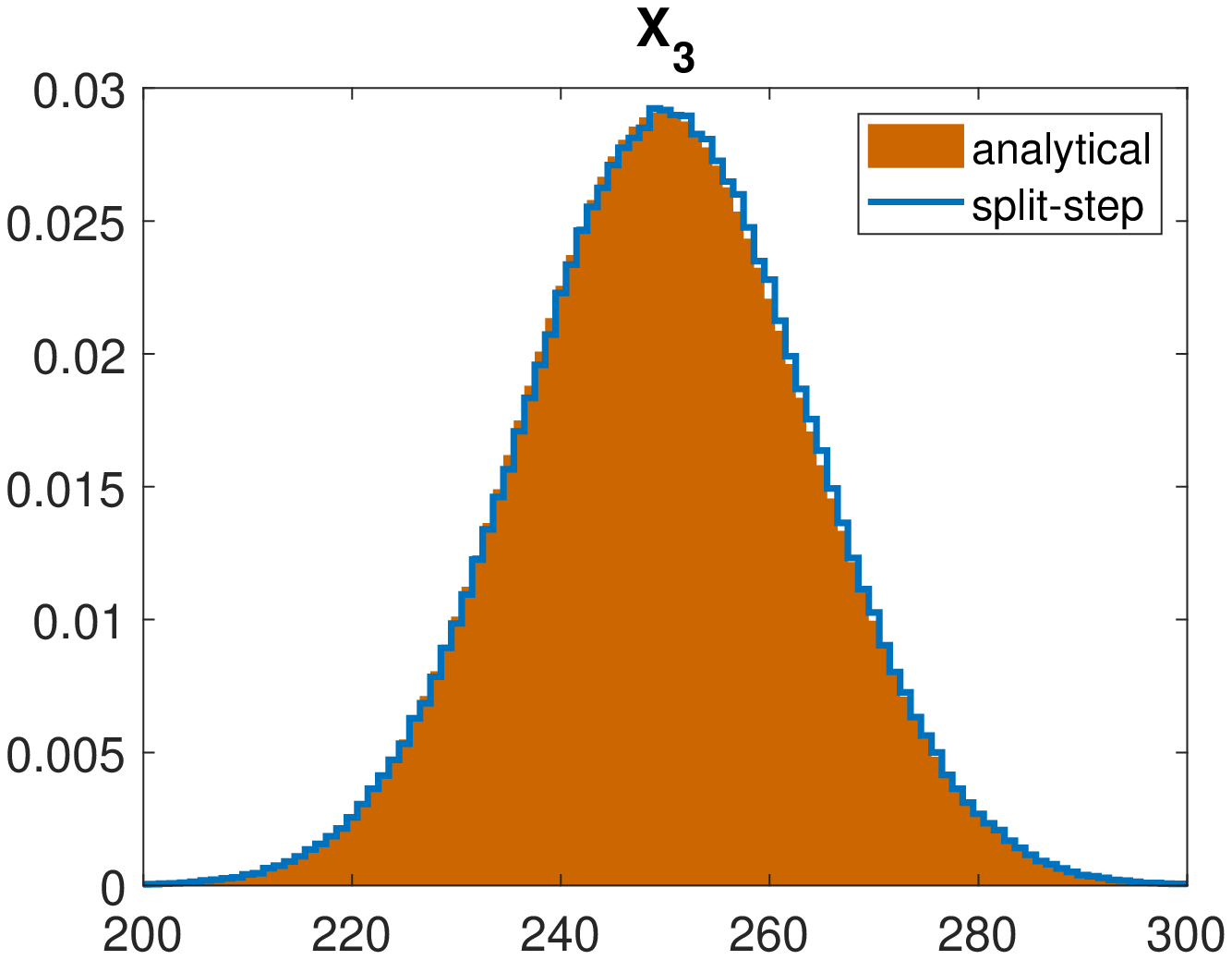}
                \caption{$x_T=1000$}
                \label{fig:ex_1_distr_xT=1000}
        \end{subfigure}
	\caption{Marginal distributions of $X_1(T)$ and $X_3(T)$ at $T=10$ for the system in \eqref{ex:1} with the reaction rates $c~=~[10^4,10^4,10^2,10^2,10^5,10^5].$}
	\label{fig:ex_1_distr}
\end{figure}

Finally, Figure \ref{fig:ex_1_distr} shows the exact and simulated marginal distributions of $X_1(T)$ and $X_3(T)$ for  different numbers of molecules $x_T$.
One can see that the split-step method \eqref{eq:spli_step} is able to accurately recover the stationary distributions for both small and large numbers of molecular species. % even at this scale.
The better accuracy for larger values of $x_T$ can be explained as an artifact of the integrality preserving procedure.

% \begin{figure}[!t]
% 	\centering
%         \begin{subfigure}[t]{1.0\textwidth}
%                 \centering
%                 \includegraphics[width=0.45\textwidth]{ex_1_xT_10_c_1e6_X1}
%                 \qquad
%                 \includegraphics[width=0.45\textwidth]{ex_1_xT_10_c_1e6_X3}
%                 \caption{$x_T=10$}
%                 \label{fig:ex_1_distr_xT=10}
%         \end{subfigure}
%         \\
%         \begin{subfigure}[t]{1.0\textwidth}
%                 \centering
%                 \includegraphics[width=0.45\textwidth]{ex_1_xT_100_c_1e6_X1}
%                 \qquad
%                 \includegraphics[width=0.45\textwidth]{ex_1_xT_100_c_1e6_X3}
%                 \caption{$x_T=100$}
%                 \label{fig:ex_1_distr_xT=100}
%         \end{subfigure}
%         \\
%         \begin{subfigure}[t]{1.0\textwidth}
%                 \centering
%                 \includegraphics[width=0.45\textwidth]{ex_1_xT_1000_c_1e6_X1}
%                 \qquad
%                 \includegraphics[width=0.45\textwidth]{ex_1_xT_1000_c_1e6_X3}
%                 \caption{$x_T=1000$}
%                 \label{fig:ex_1_distr_xT=1000}
%         \end{subfigure}
% 	\caption{Marginal distributions of $X_1(T)$ and $X_3(T)$ at $T=100$ for the system in \eqref{ex:1} with the reaction rates $c~=~10^6 \cdot [1,1,1,1,1,1].$}
% 	\label{fig:ex_1_distr}
% \end{figure}

\subsection{Example 2. (Nonlinear system)}

Consider the following stiff 3-species 6-reaction system \cite{Sotiropoulos2008}

\begin{align}\label{ex:2}
	\nonumber
	S_1 + S_2 \underset{c_2}{\stackrel{c_1}{\rightleftharpoons}} S_3,
	\\
	S_1 + S_3 \underset{c_4}{\stackrel{c_3}{\rightleftharpoons}} S_2,
	\\ \nonumber
	S_2 + S_3 \underset{c_6}{\stackrel{c_5}{\rightleftharpoons}} S_1
\end{align}
with the rate constants $c_1 = c_2=10^3$, $c_3=10^{-5}$, $c_4=10$, $c_5=1$ and $c_6=10^6$.
The stoichiometric matrix of this reaction network is given by

\begin{align*}
	\nu = 
	\begin{pmatrix}
		-1 & 1 & -1 & 1 & 1 & -1 \\[0.5em]
		-1 & 1 & 1 & -1 & -1 & 1 \\[0.5em]
		1 & -1 & -1 & 1 & -1 & 1 \\
	\end{pmatrix}
\end{align*}
and the propensities of the reaction channels are, respectively,

\begin{alignat*}{2}
	&\alpha_1(X(t)) = c_1 X_1(t) X_2(t), \quad & \alpha_2(X(t)) = c_2 X_3(t), & \\
	&\alpha_3(X(t)) = c_3 X_1(t) X_3(t), \quad & \alpha_4(X(t)) = c_4 X_2(t), & \\
	&\alpha_5(X(t)) = c_5 X_2(t) X_3(t), \quad & \alpha_6(X(t)) = c_6 X_1(t). & 
\end{alignat*}

The system in \eqref{ex:2} does not have known analytical solutions.
Instead, the reference solution was obtained using the stochastic simulation algorithm with $10^4$ samples.
We used the initial condition $X(0) = [10^3, 10^3, 10^6]^T$ which corresponds to the true stationary mean of this particular system and performed all simulations over the time interval $[0,0.01]$ with the time step $\tau=10^{-3}$.

\begin{figure}[!t]
	\centering
        \begin{subfigure}[t]{1.0\textwidth}
                \centering
                \includegraphics[width=0.45\textwidth]{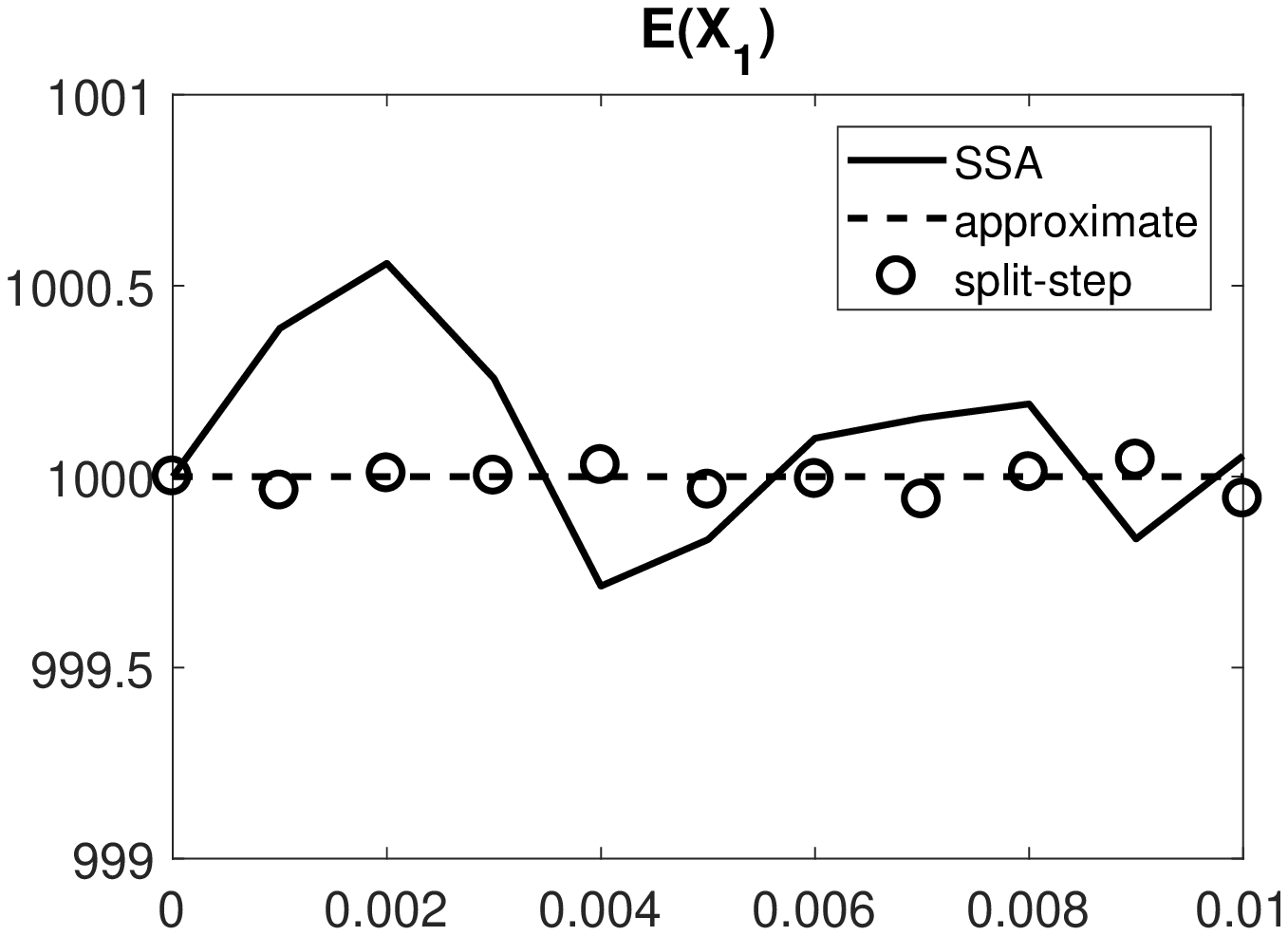}
                \qquad
                \includegraphics[width=0.45\textwidth]{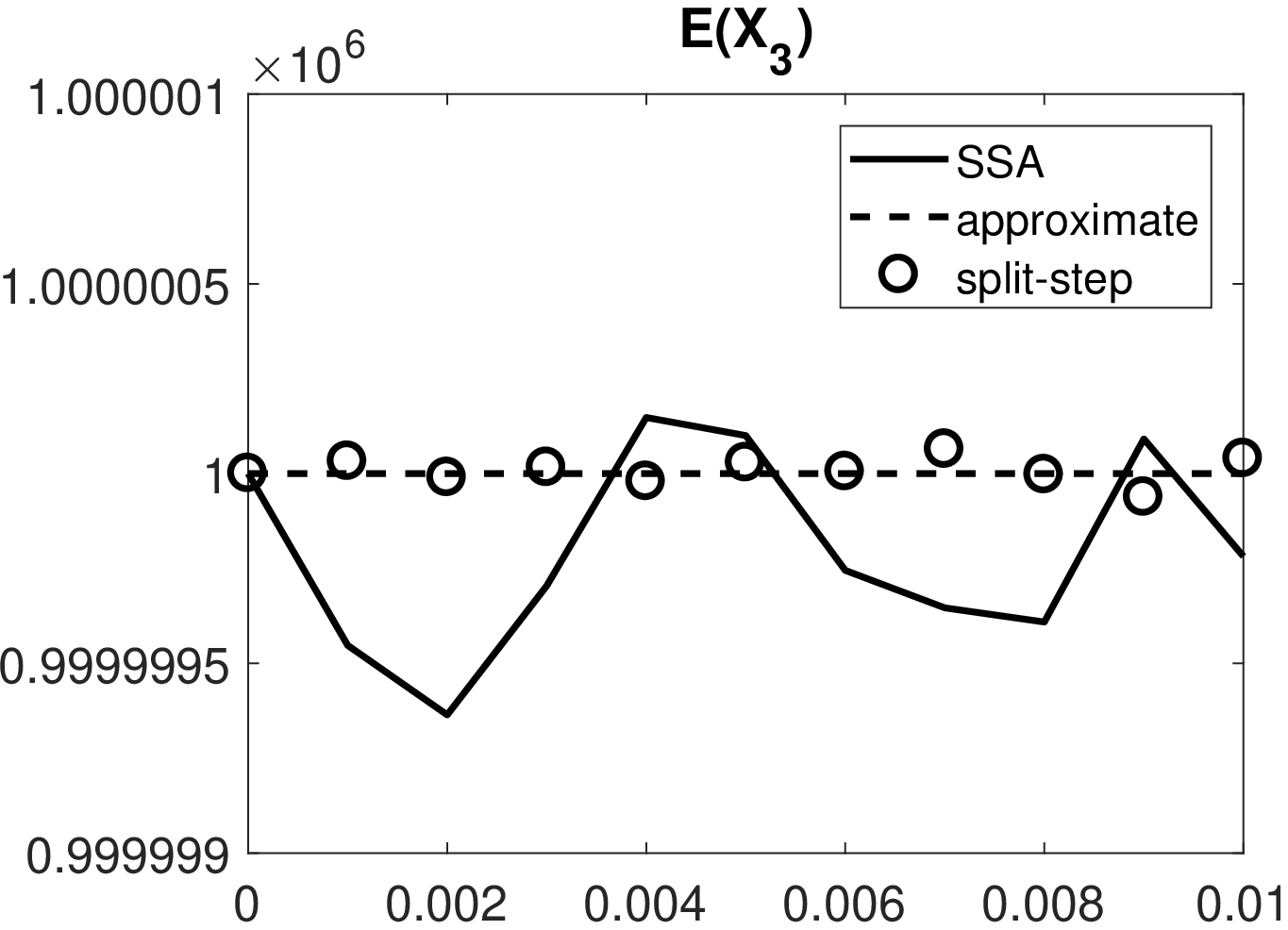}
        \end{subfigure}
        \\
        \begin{subfigure}[t]{1.0\textwidth}
                \centering
                \includegraphics[width=0.45\textwidth]{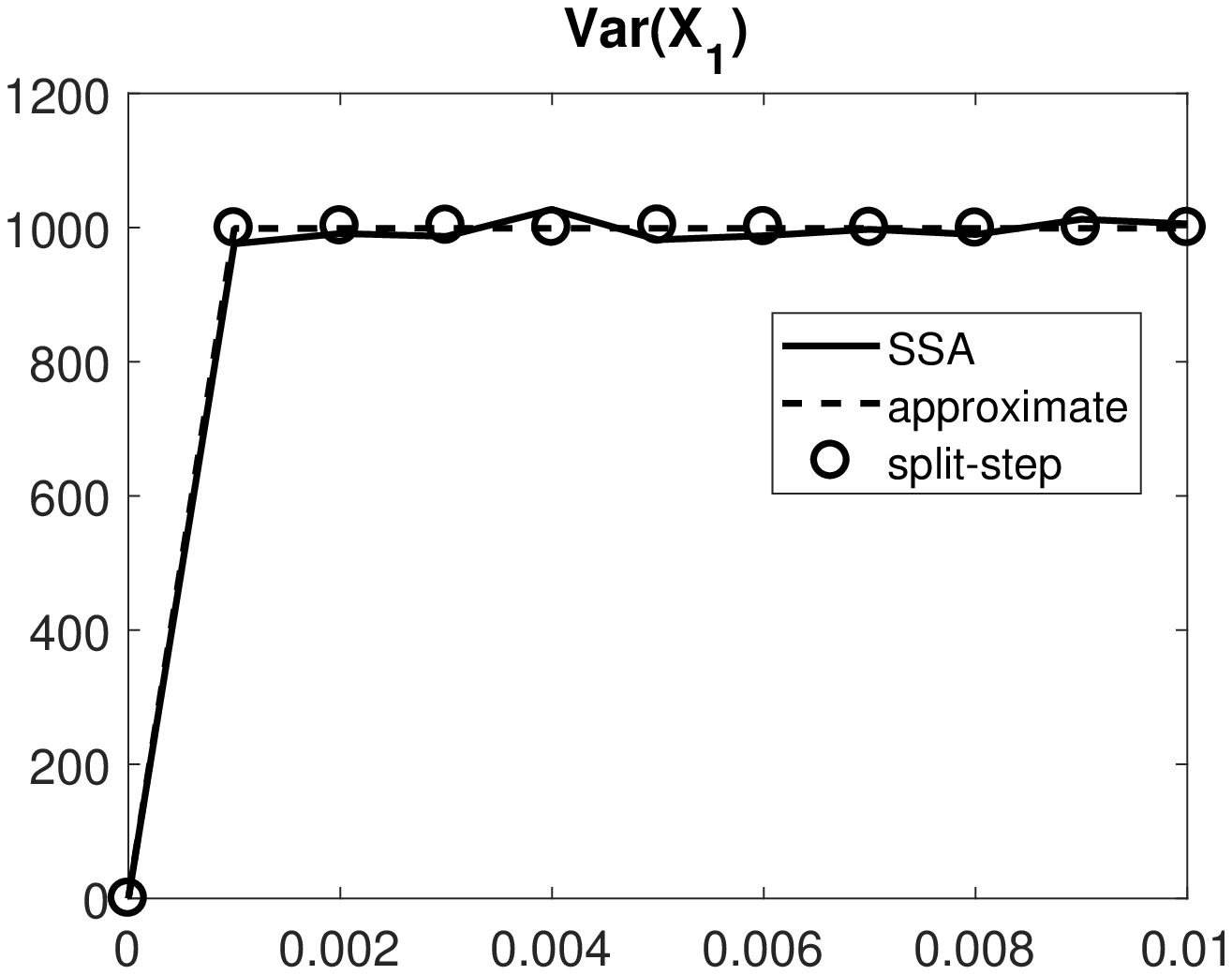}
                \qquad
                \includegraphics[width=0.45\textwidth]{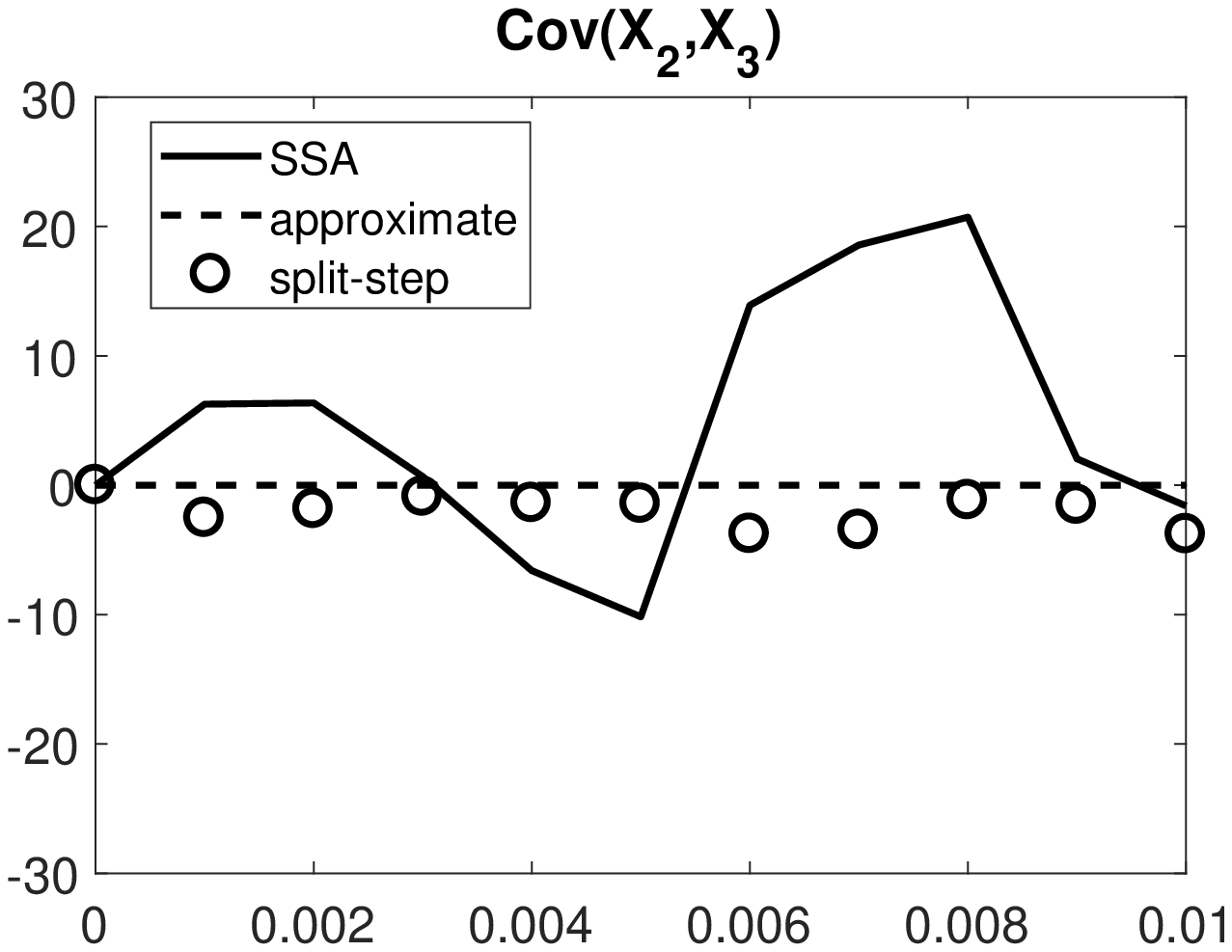}
        \end{subfigure}
%         \\
%         \begin{subfigure}[t]{1.0\textwidth}
%                 \centering
%                 \includegraphics[width=0.45\textwidth]{ex_2_cov13}
%                 \qquad
%                 \includegraphics[width=0.45\textwidth]{ex_2_cov23}
%         \end{subfigure}
	\caption{Evolution in time of the mean and covariance of the system in \eqref{ex:2}.}
	\label{fig:ex_2_mean_cov}
\end{figure}

Figure \ref{fig:ex_2_mean_cov} depicts the evolution in time of the selected elements of the mean vector and the covariance matrix.
The approximate mean and covariance were evaluated using the linearization in \eqref{eq:linearized_prop}.
It is seen that both the linearization \eqref{eq:linearized_prop} and the corresponding fitted approximation \eqref{eq:split_step_mean}-\eqref{eq:split_step_var} produce an adequate description of the true stationary nonlinear dynamics.

\begin{figure}[!t]
	\centering
        \begin{subfigure}[t]{1.0\textwidth}
                \centering
                \includegraphics[width=0.45\textwidth]{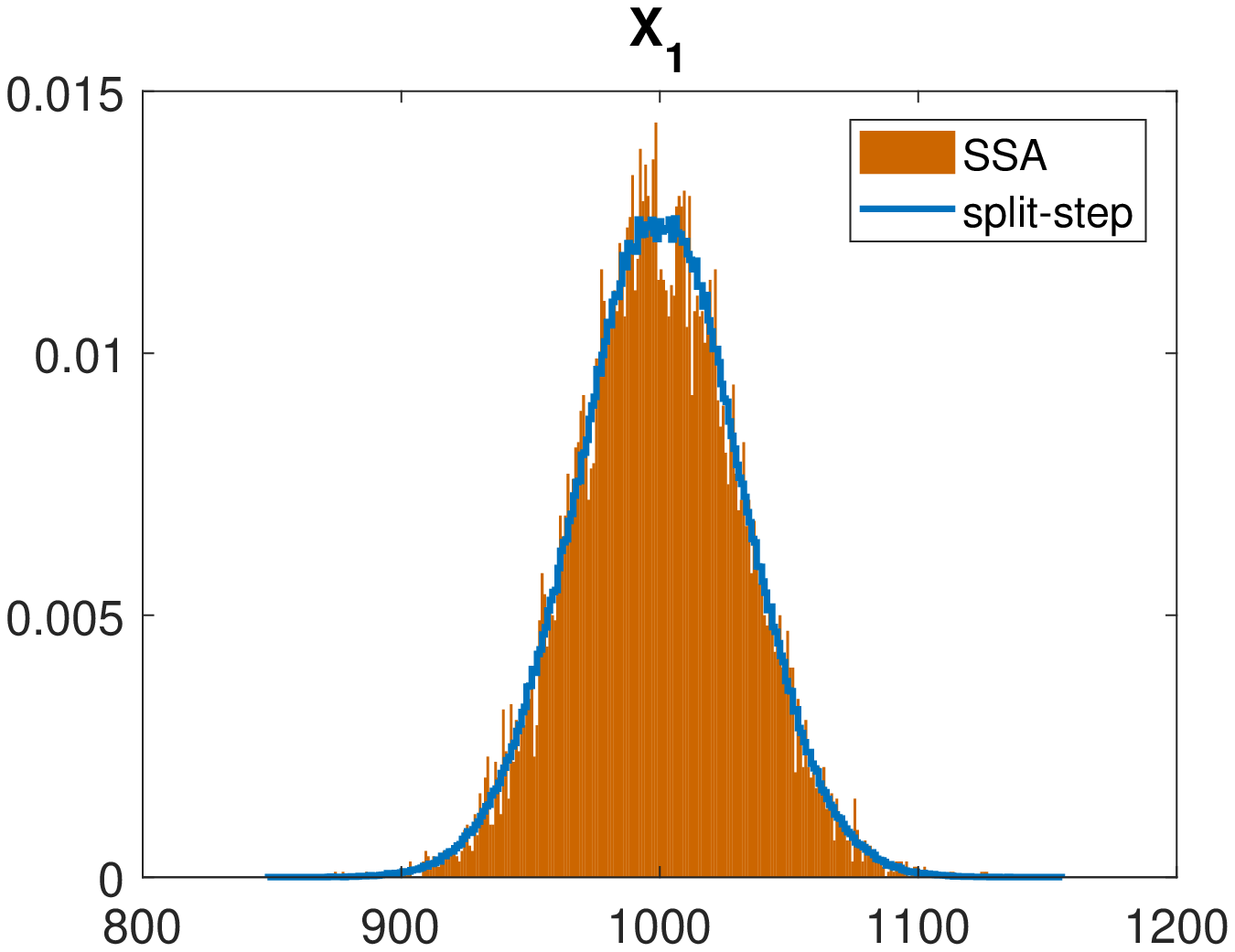}
                \qquad
                \includegraphics[width=0.45\textwidth]{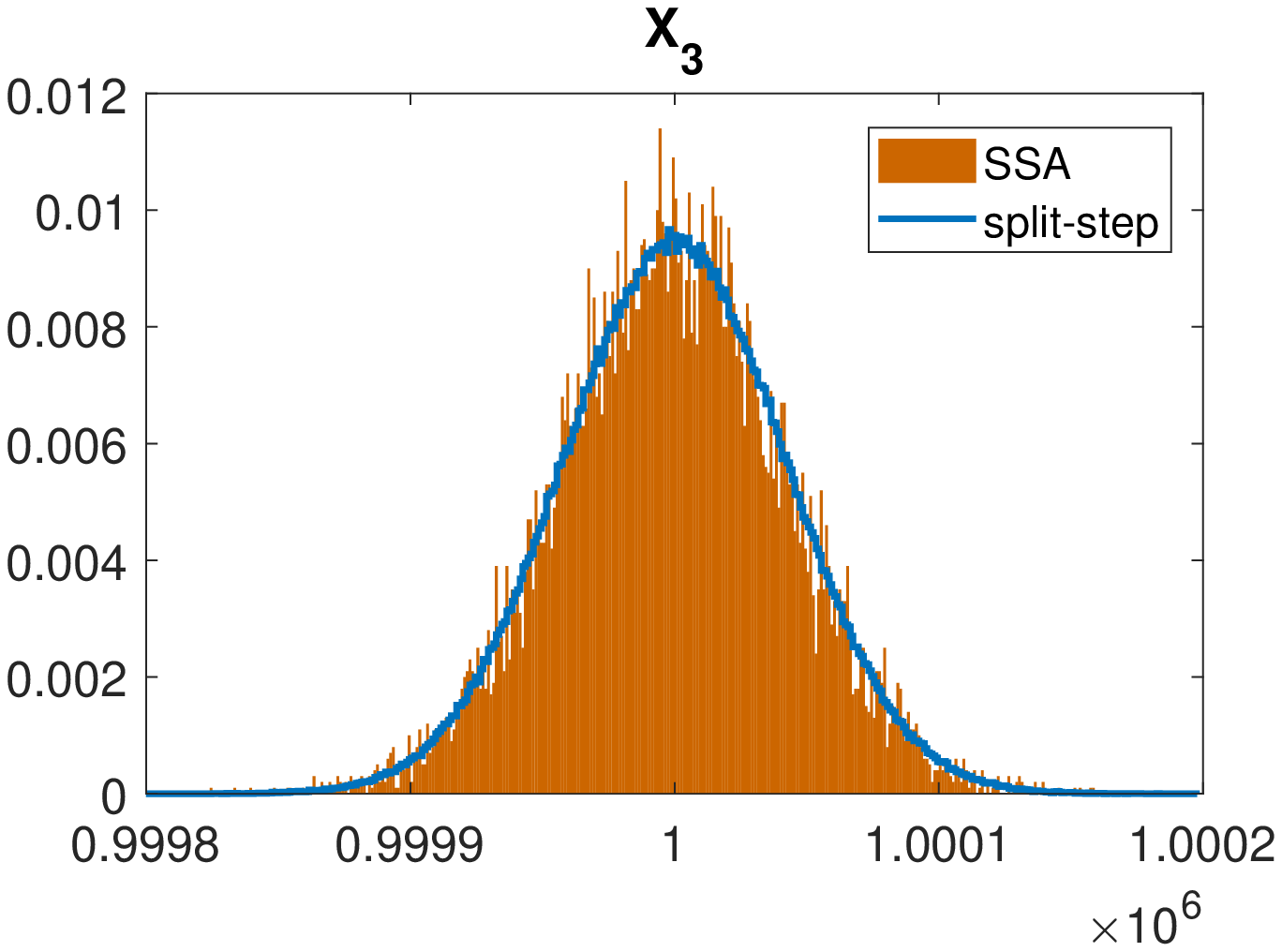}
        \end{subfigure}
        \\
        \begin{subfigure}[t]{1.0\textwidth}
                \centering
                \includegraphics[width=0.45\textwidth]{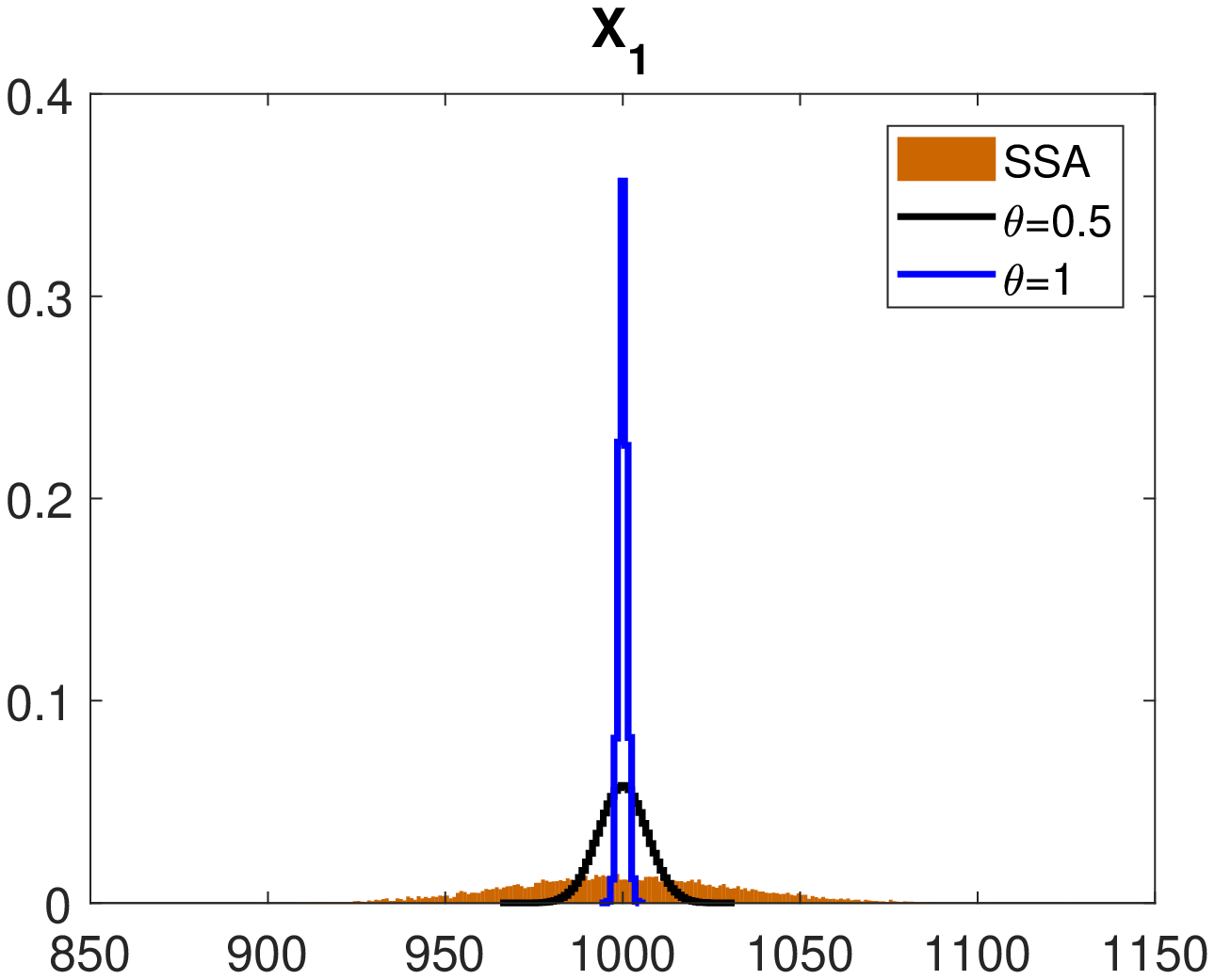}
                \qquad
                \includegraphics[width=0.45\textwidth]{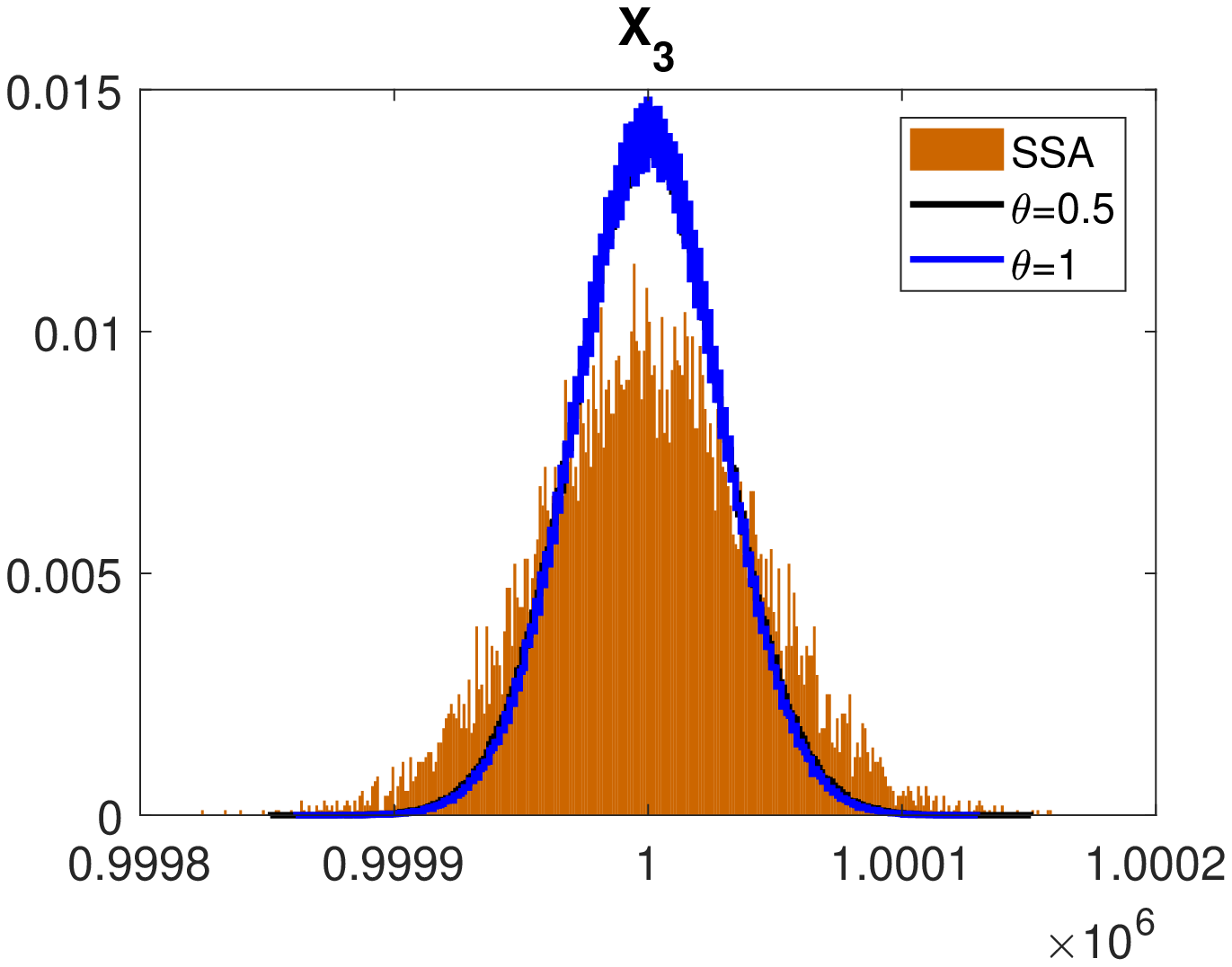}
        \end{subfigure}
	\caption{Marginal distributions of $X_1(T)$ and $X_3(T)$ at $T=0.01$ for the system in \eqref{ex:2}.}
	\label{fig:ex_2_distr}
\end{figure}

Additionally, Figure \ref{fig:ex_2_distr} provides the comparison of the marginal distributions obtained with the split-step and theta schemes.
It is evident that both implicit and trapezoidal methods produce severely overdamped solutions resulting in very atomic distributions.
At the same time, the proposed split-step integrator generates results which compare well to the reference solution as desired.
% The proposed integrator gives result which compare very well with the reference solution as desired.

\section{Conclusion}

\change{
In this paper, we considered the problem of accurate integration of the fast and stable reactions using time discretizations much exceeding their scale.
For this purpose, we proposed the new splitting heuristic for the tau-leaping method and devised the algorithm for the estimation of its parameters.
We showed that the proposed integrator is stable and is capable to accurately sample from the stationary distributions of the fast molecular species at least in the linear case.
%In this paper, we proposed the new splitting heuristic for the tau-leaping methods and developed the algorithm for the estimation of its optimal parameters.
%The considered split-step technique has the advantage of being L-stable in the mean while accurately resolving the covariance structure of stationary solutions.
Classical tau-leaping approximations do not possess this property and behave poorly as stiff integrators.
%This feature is particularly useful for the integration of stiff stochastic systems since it allows to step over fast reactions and capture their influence on slow variables.

% The parameters of the method were calibrated using linear moment analysis but the numerical results indicate that the method works well for both linear and nonlinear systems with well defined first two moments.
% It was also observed that the mean and the covariance can be accurately resolved for a large range of molecular counts but the accuracy in the distribution deteriorates for the systems with small number of species.
% Finally, the proposed algorithm does not assume any particular type of equations and dynamics but shows the best results for the systems in equilibrium and can behave poorly for simulating fast transients.

%It is worth noting that the proposed split-step tau-leaping method can be used as a self-contained stiff integrator but also as a low-fidelity solver in the multilevel setting \cite{BenHammouda2017,Peherstorfer2016,Giles2008}.
}

\appendix
\section{Stability of the theta method}
\label{sec:theta_stab}

Let $Y_n$ denote the numerical approximation of $X_1(t_n)$ and consider the theta method \eqref{eq:theta} applied to the test system \eqref{eq:test}
\begin{align*}
	Y_{n+1}
	&= Y_n + \Big( c_2 (x_T - Y_{n+1}) - c_1 Y_{n+1} - c_2 (x_T - Y_{n}) + c_1 Y_{n} \Big) \theta\tau 
	\\
	&+\mathcal{P}\Big( c_2 (x_T - Y_{n}) \tau \Big) - \mathcal{P}\Big( c_1 Y_{n} \tau \Big)
	\\
	&= Y_n + \frac{1}{1+(c_1+c_2)\theta\tau} \left[ \mathcal{P}\Big( c_2 (x_T - Y_{n}) \tau \Big) - \mathcal{P}\Big( c_1 Y_{n} \tau \Big) \right].
\end{align*}	
The mean and the variance of the solution obtained with the theta method can be easily calculated as
\begin{align*}
	\E{Y_{n+1}} &= \E{\E{Y_{n+1} | Y_n}}
	\\
	&= \E{Y_n + \frac{c_2 (x_T - Y_{n}) \tau - c_1 Y_{n} \tau}{1+(c_1+c_2)\theta\tau} }
	= \frac{1 - (c_1+c_2)(1-\theta)\tau}{1+(c_1+c_2)\theta\tau} \E{Y_n} + \frac{c_2 x_T \tau}{1+(c_1+c_2)\theta\tau}
\end{align*}
and
\begin{align*}
	\Var{Y_{n+1}} &= \Var{\E{Y_{n+1} | Y_n}} + \E{\Var{Y_{n+1} | Y_n}}
	\\[1em]
	&= \Var{ Y_n + \frac{c_2 (x_T - Y_{n}) \tau - c_1 Y_{n} \tau}{1+(c_1+c_2)\theta\tau} }
	+ \E{ \frac{c_1 Y_{n} \tau  +  c_2 (x_T - Y_{n}) \tau}{(1+(c_1+c_2)\theta\tau)^2} }
	\\[1em]
	&= \left( \frac{1 - (c_1+c_2)(1-\theta)\tau}{1+(c_1+c_2)\theta\tau} \right)^2 \Var{Y_n} 
	+  \frac{(c_1-c_2)\tau}{(1+(c_1+c_2)\theta\tau)^2} \E{Y_n} 
	+  \frac{c_2 x_T \tau}{(1+(c_1+c_2)\theta\tau)^2}.
\end{align*}
These recurrences have the solutions
\begin{align*}
	\E{Y_{n}} &= 
	%A^n \E{Y_0} + \left( \sum_{i=0}^{n-1} A^i \right) B = 
	A^n \E{Y_0} + \frac{1-A^n}{1-A} B,
	\\[1em]
	\Var{Y_{n}} &= A^{2n} \Var{Y_0} + \frac{C}{A} \frac{1-A^n}{1-A} \left( \E{Y_{n}} - \frac{B}{1-A} \right) + \frac{1-A^{2n}}{1-A^2} \left( \frac{B C}{1-A} + D \right),
\end{align*}
where the coefficients $A$, $B$, $C$ and $D$ are given by
\begin{alignat*}{2}
	& A=\frac{1 - (c_1+c_2)(1-\theta)\tau}{1+(c_1+c_2)\theta\tau}, \qquad && B=\frac{c_2 x_T \tau}{1+(c_1+c_2)\theta\tau}, 
	\\
	&C=\frac{(c_1-c_2)\tau}{(1+(c_1+c_2)\theta\tau)^2}, \qquad && D=\frac{c_2 x_T \tau}{(1+(c_1+c_2)\theta\tau)^2}.
\end{alignat*}

The following condition ensures the global stability of the mean and the variance of the numerical solution obtained with the theta method \eqref{eq:theta}
\begin{align*}
	|A| = \left| \frac{1 - (c_1+c_2)(1-\theta)\tau}{1+(c_1+c_2)\theta\tau} \right| < 1.
\end{align*}
%Assuming that $|A|<1$ and letting $n\to\infty$, we get the mean and the variance of the stationary distribution generated by the theta method
%\begin{align}
%	\E{Y_{\infty}} &= \frac{B}{1-A} 
%	%= \frac{c_2}{c_1+c_2}x_T 
%	= \E{X^{*}},
%	\\ \nonumber
%	\Var{Y_{\infty}} &= \frac{1}{1-A^2} \left( \frac{B C}{1-A} + D \right) 
%	%= \frac{2}{2+(c_1+c_2)(2 \theta-1)\tau} \frac{c_1 c_2 x_T}{(c_1+c_2)^2} 
%	= \frac{2}{2+(c_1+c_2)(2 \theta-1)\tau} \Var{X^{*}}.
%\end{align}

\section{Stability of the split-step method}
\label{sec:split_step_stab}
Let $Y_n$ denote the numerical approximation of $X_1(t_n)$ and consider the split-step method \eqref{eq:mod_split_step} applied to the test system \eqref{eq:test}
\begin{align*}
	\hat{Y}_{n}
	&= \frac{1-(1-\eta_1)(c_1+c_2)(1-\theta)\tau}{1+\eta_1(c_1+c_2)(1-\theta)\tau} Y_n + \frac{ c_2 x_T (1-\theta) \tau}{1+\eta_1(c_1+c_2)(1-\theta)\tau} ,
	\\
	\tilde{Y}_{n}
	&= (1+\tau(c_1+c_2))\hat{Y}_{n} + \mathcal{P}\left(c_2(x_T-\hat{Y}_{n}) \tau \right) - \mathcal{P}\left(c_1\hat{Y}_{n} \tau \right) - c_2 \tau x_T,
	\\
	Y_{n+1}
	&= \frac{1-(1-\eta_2)(c_1+c_2)\theta\tau}{1+\eta_2(c_1+c_2)\theta\tau} \tilde{Y}_{n} + \frac{c_2 x_T \theta\tau}{1+\eta_2(c_1+c_2)\theta\tau}.
\end{align*}
Equations for the mean and the variance take the form
\begin{align*}
	\E{Y_{n+1}}
	&=   \left(\frac{1-(1-\eta_1)(c_1+c_2)(1-\theta)\tau}{1+\eta_1(c_1+c_2)(1-\theta)\tau} \right) \left( \frac{1-(1-\eta_2)(c_1+c_2)\theta\tau}{1+\eta_2(c_1+c_2)\theta\tau} \right)  \E{Y_n} 
	\\[1em]
	&+\frac{1 - (1-\eta_1-\eta_2)(c_1+c_2)\theta(1-\theta)\tau}{\Big(1+\eta_1(c_1+c_2)(1-\theta)\tau\Big) \Big(1+\eta_2(c_1+c_2)\theta\tau\Big)} c_2 x_T \tau,
	\\[1em]
	\Var{Y_{n+1}}
	&=  \left( \frac{1-(1-\eta_1)(c_1+c_2)(1-\theta)\tau}{1+\eta_1(c_1+c_2)(1-\theta)\tau} \right)^2 \left( \frac{1-(1-\eta_2)(c_1+c_2)\theta\tau}{1+\eta_2(c_1+c_2)\theta\tau} \right)^2  \Var{Y_n} 
	\\[1em]
	&+ \left(\frac{1-(1-\eta_1)(c_1+c_2)(1-\theta)\tau}{1+\eta_1(c_1+c_2)(1-\theta)\tau} \right) \left( \frac{1-(1-\eta_2)(c_1+c_2)\theta\tau}{1+\eta_2(c_1+c_2)\theta\tau} \right)^2  (c_1-c_2) \tau  \E{Y_n}
	\\[1em]
	&+ \left( \frac{1-(1-\eta_2)(c_1+c_2)\theta\tau}{1+\eta_2(c_1+c_2)\theta\tau} \right)^2 \left( 1 + \frac{ (c_1-c_2) (1-\theta) \tau}{1+\eta_1(c_1+c_2)(1-\theta)\tau} \right) c_2 x_T \tau
\end{align*}	
and the global stability condition reads as
\begin{align*}
	\left|  \left(\frac{1-(1-\eta_1)(c_1+c_2)(1-\theta)\tau}{1+\eta_1(c_1+c_2)(1-\theta)\tau} \right) \left( \frac{1-(1-\eta_2)(c_1+c_2)\theta\tau}{1+\eta_2(c_1+c_2)\theta\tau} \right)  \right| < 1.
\end{align*}	

\section*{References}

\bibliography{references}

\end{document}